\documentclass[final,leqno,onefignum]{siamltex704}

\usepackage{fancyhdr}
\usepackage{mathrsfs}
\usepackage{latexsym}
\usepackage{amsfonts}
\usepackage[usenames]{color}
\usepackage{graphicx}
\usepackage{amsmath}
\usepackage{enumitem}
\newtheorem{remark}{Remark}[section]
\newtheorem{assumption}{Assumption}[section]

\newcommand{\be}[1]{\begin{equation} \label{#1}}
\newcommand{\ee}{\end{equation}}
\newcommand{\pdpd}[2]{\frac{\partial #1}{\partial #2}}

\newcommand{\rmd}{\mathrm{d}}
\newcommand{\ii}{\mathrm{i}}
\newcommand{\R}{\mathbb{R}}

\renewcommand{\Re}{\operatorname{Re}}
\newcommand{\half}{\mbox{$\frac{1}{2}$}}

\def \p {\partial}
\def \supp  {\operatorname{supp}}
\def \op {\operatorname{Op}}

\def \p {\partial}

\def \vphi {\varphi}

\begin{document}

\title{Semiclassical analysis of elastic surface waves}
\author{Maarten V. de Hoop \thanks{Simons Chair in Computational and
    Applied Mathematics and Earth Science, Rice University, Houston,
    TX, 77005, USA (mdehoop@rice.edu)}
\and Alexei Iantchenko \thanks{Department of Materials Science and Applied Mathematics, Faculty of Technology and Society, Malm\"{o} University, SE-205 06 Malm\"{o}, Sweden (ai@mah.se)}
\and Gen Nakamura \thanks{Department of Mathematics, Hokkaido University,
 Sapporo 060-0810, Japan (nakamuragenn@gmail.com)}
\and Jian Zhai \thanks{Institute for Advanced Study, The Hong Kong University of Science and Technology, Kowloon, Hong Kong, China
  (jian.zhai@outlook.com)}}

\maketitle

\pagestyle{myheadings}
\thispagestyle{plain}
\markboth{DE HOOP, IANTCHENKO, NAKAMURA, and ZHAI}{Semiclassical
  analysis of elastic surface waves}

\begin{abstract}
In this paper, we present a semiclassical description of surface waves
or modes in an elastic medium near a boundary, in spatial dimension
three. The medium is assumed to be essentially stratified near the
boundary at some scale comparable to the wave length. Such a medium
can also be thought of as a surficial layer (which can be thick)
overlying a half space. The analysis is based on the work of Colin de
Verdi\`ere \cite{CDV2} on acoustic surface waves. The description is
geometric in the boundary and locally spectral ``beneath''
it. Effective Hamiltonians of surface waves correspond with
eigenvalues of ordinary differential operators, which, to leading
order, define their phase velocities. Using these Hamiltonians, we
obtain pseudodifferential surface wave equations. We then construct a
parametrix. Finally, we discuss Weyl's formulas for counting surface
modes, and the decoupling into two classes of surface waves, that is,
Rayleigh and Love waves, under appropriate symmetry conditions.
\end{abstract}

\section{Introduction}
\label{sec:1}

We carry out a semiclassical analysis of surface waves in a medium
which is stratified near its boundary -- with topography -- at some
scale comparable to the wave length. We discuss how the (dispersive)
propagation of such waves is governed by effective Hamiltonians on the
boundary and show that the system is displayed by a space-adiabatic
behavior.

The Hamiltonians are non-homogeneous principal symbols of some
pseudodifferential operators. Each Hamiltonian is identified with an
eigenvalue in the point spectrum of a locally Schr\"{o}dinger-like
operator in dimension one on the one hand, and generates a flow
identified with surface-wave bicharacteristics in the two-dimensional
boundary on the other hand. The eigenvalues exist under certain
assumptions reflecting that wave speeds near the boundary are smaller
than in the deep interior. This assumption is naturally satisfied by
the structure of Earth's crust and mantle (see, for example, Shearer
\cite{shearer}). The dispersive nature of surface waves is manifested
by the non-homogeneity of the Hamiltonians.

The spectra of the mentioned Schr\"{o}dinger-like operators consist of
point and essential spectra. The surface waves are identified with the
point spectra while the essential spectra correspond with propagating
body waves. We note, here, that the point and essential spectra for
the Schr\"{o}dinger-like operators may overlap.

Our analysis applies to the study of surface waves in Earth's ``near''
surface in the scaling regime mentioned above. The existence of such
waves, that is, propagating wave solutions which decay exponentially
away from the boundary of a homogeneous (elastic) half-space was first
noted by Rayleigh \cite{Rayleigh}. Rayleigh and (``transverse'') Love
waves can be identified with Earth's free oscillation triples ${}_n
S_l$ and ${}_n T_l$ with $n \ll l/4$ assuming spherical symmetry. Love
\cite{Love} was the first to argue that surface-wave dispersion is
responsible for the oscillatory character of the main shock of an
earthquake tremor, following the ``primary'' and ``secondary''
arrivals.

Our analysis is motivated by the (asymptotic) JWKB theory of surface
waves developed in seismology by Woodhouse \cite{Woodhouse}, Babich,
Chichachev and Yanoskaya \cite{Babich} and others. Tromp and Dahlen
\cite{Tromp} cast this theory in the framework of a ``slow''
variational principle. The theory is also used in ocean acoustics
\cite{Brekhovskikh} and is referred to as adiabatic mode theory. An
early study of the propagation of waves in smoothly varying waveguides
can be found in Bretherton \cite{Bretherton}. Nomofilov \cite{Nomo}
obtained the form of WKB solutions for Rayleigh waves in
inhomogeneous, anisotropic elastic media using assumptions appearing
in Proposition~\ref{R_speed} in the main text. Many aspects of the
propagation of surface waves in laterally inhomogeneous elastic media
are discussed in the book of Malischewsky \cite{Malischewsky}. Here,
we develop a comprehensive semiclassical analysis of elastic surface
waves, generated by interior (point) sources, with the corresponding
estimates. This semiclassical framework was first formulated by Colin
de Verdi\`ere \cite{CDV2} to describe surface waves in acoustics.

The scattering of surface waves by structures, away from the mentioned
scaling regime, has been extensively studied in the seismology
literature. This scattering can be described using a basis of local
surface wave modes that depend only on the ``local'' structure of the
medium, for example, with invariant embedding; see, for example, Odom
\cite{Odom}. Odom used a layer of variable thickness over a
homogeneous half space to account for the interaction between surface
waves and body waves by the topography of internal
interfaces.

The outline of this paper is as follows. In Section~\ref{sec:2}, we
carry out the semiclassical construction of general surface wave
parametrices. In the process, we introduce locally
Schr\"{o}dinger-like operators in the boundary normal coordinate and
their eigenvalues signifying effective Hamiltonians in the boundary
(tangential) coordinates describing surface-wave propagation. In
Section~\ref{sec:3}, we characterize the spectra of the relevant
Schr\"{o}dinger-like operators. That is, we study their discrete and
essential spectra. In Section~\ref{sec:4}, we consider a special class
of surface modes associated with exponentially decaying
eigenfunctions. The existence of such modes is determined by a
generalized Barnett-Lothe condition. In Section~\ref{sec:5}, we review
conditions on the symmetry, already considered by Anderson
\cite{Anderson}, restricting the anisotropy to transverse isotropy
with the axis of symmetry aligned with the normal to the boundary,
allowing the decoupling of surface waves into Rayleigh and Love
waves. In Section~\ref{sec:6}, we establish Weyl's laws first in the
isotropic (separately for Rayleigh and Love waves) and then in the
anisotropic case. Finally, in Section~\ref{sec:7}, we relate the
surface waves to normal modes viewing the analysis locally on conic
regions, or more generally on Riemannian manifolds with a half
cylinder structure. We give explicit formulas for the special case of
a radial manifold.

\section{Semiclassical construction of surface-wave parametrices}
\label{sec:2}

We consider the linear elastic wave equation in $\mathbb{R}^3$,
\begin{equation}\label{Elastic1}
 \frac{\partial^2 u}{\partial t^2}=\mathrm{div}~\frac{\sigma(u)}{\rho},
\end{equation}
where $u$ is the displacement vector, and $\sigma(u)$ is the stress tensor given by Hooke's law
\begin{equation}\label{Hooke}\sigma(u)=\mathbf{C}\varepsilon (u),\end{equation}
and $\varepsilon(u)$ denotes the strain tensor; $\mathbf{C}$ is the fourth-order stiffness tensor with components $c_{ijkl}$, and $\rho$ is the density of mass.
The componentwise expression of $(\ref{Hooke})$ is given by
\[\sigma_{ij}(u)=\sum_{k,l=1}^3c_{ijkl}\varepsilon_{kl}(u),~\varepsilon_{kl}(u)=\frac{1}{2}(\partial_ku_l+\partial_lu_k).\]
Equation (\ref{Elastic1}) differs from the usual system given by
\[
   \rho\frac{\partial^2 u}{\partial t^2} =\mathrm{div}~\sigma(u)
\]
in case $\rho$ is not a constant. However, the difference is in the
lower terms. These can be accounted for considering (\ref{Elastic1}).

We study the elastic wave equation $(\ref{Elastic1})$ in the half space $X = \mathbb{R}^2
\times (-\infty,0]$, with coordinates,
\[
(x,z) ,\quad x = (x_1,x_2) \in \mathbb{R}^2 ,\
z \in \mathbb{R}^{-} = (-\infty,0] .
\]
We consider solutions, $u = (u_1,u_2,u_3)$, satisfying the Neumann
boundary condition at $\partial X = \{z=0\}$,
\begin{equation}\label{elaswaeq}
\begin{split}
   \partial^2_t u_i + M_{il} u_l &= 0 ,\\
   u(t=0,x,z) &= 0,~~\partial_tu(t=0,x,z)=h(x,z) ,\\
   \frac{c_{i3kl}}{\rho}\partial_k u_l(t,x,z=0) &= 0 ,
\end{split}
\end{equation}
where
\begin{multline*}
M_{il} =
-\frac{\partial}{\partial z}\frac{c_{i33l}(x,z)}{\rho(x,z)}
\frac{\partial}{\partial z}     - \sum_{j,k=1}^{2}
\frac{c_{ijkl}(x,z)}{\rho(x,z)} \frac{\partial}{\partial x_j}
\frac{\partial}{\partial x_k}
- \sum_{j=1}^{2}
\frac{\partial}{\partial x_j}\frac{c_{ij3l}(x,z)}{\rho(x,z)}
\frac{\partial}{\partial z}
\\
- \sum_{k=1}^{2}
\frac{c_{i3kl}(x,z)}{\rho(x,z)} \frac{\partial}{\partial z}
\frac{\partial}{\partial x_k}
-\sum_{k=1}^{2}
\left( \frac{\partial}{\partial z} \frac{c_{i3kl}(x,z)}{\rho(x,z)} \right)
\frac{\partial}{\partial x_k}
\\
- \sum_{j,k=1}^{2}
\left( \frac{\partial}{\partial x_j} \frac{c_{ijkl}(x,z)}{\rho(x,z)} \right)
\frac{\partial}{\partial x_k}.
\hspace*{2.5cm}
\end{multline*}
In the above, we assume that $\rho(x,z) \in C^\infty(X)$ with
$\rho(x,z) \ge \rho_0 > 0$ for some $\rho_0 > 0$ and $c_{ijkl}(x,z)
\in C^\infty(X)$ satisfies the following symmetries and strong
convexity condition: \\[0.25cm]
(symmetry) $c_{ijkl}=c_{jikl}=c_{klij}$ for any $i,j,k,l$; \\[0.25cm]
(strong convexity) there exists $\delta>0$ such that for any nonzero
$3\times 3$ real-valued symmetric matrix $(\varepsilon_{ij})$,
\[
   \sum_{i,j,k,l=1}^3\frac{c_{ijkl}}{\rho}\varepsilon_{ij}\varepsilon_{kl}\geq\delta\sum_{i,j=1}^3\varepsilon_{ij}^2.
\]
We note that these are physically very natural assumptions. We invoke,
additionally,

\begin{assumption}\label{assu_tensor}
The stiffness tensor and density obey the following scaling,
\[
   \frac{c_{ijkl}}{\rho}(x,z)
     = C_{ijkl}\left(x,\frac{z}{\epsilon}\right) ,
                            ~~\epsilon\in(0,\epsilon_0] ;
\]
\[
   C_{ijkl}(x,Z)=C_{ijkl}(x,Z_I)=C_{ijkl}^I
                            = \text{constant, for }Z \leq Z_I<0
\]
and, for any $\hat{\xi} \in S^2$,
\[
   \inf_{Z \leq 0}v_L(x,\hat{\xi},Z)<v_L(x,\hat{\xi},Z_I),
\]
\end{assumption}

\noindent
where $v_L$ is the so-called limiting velocity which will be defined
in Section~\ref{sec:3}.

\subsection{Schr\"{o}dinger-like operators}

Under Assumption \ref{assu_tensor}, we make the following change of
variables,
\[
   u(t,x,z) = v\left(t,x,\frac{z}{\epsilon}\right) ;
\]
upon introducing $Z = \frac{z}{\epsilon}$, the elastic wave equation
in (\ref{elaswaeq}) takes the form
\begin{equation}\label{semi1}
   [\epsilon^2 \partial_t^2 + \hat{H}] v = 0 ,
\end{equation}
where
\begin{multline}\label{Full}
\hat{H}_{il} =
-\frac{\partial}{\partial Z} C_{i33l}(x,Z)
\frac{\partial}{\partial Z}
\\
- \epsilon \sum_{j=1}^{2}
C_{ij3l}(x,Z)\frac{\partial}{\partial x_j}
\frac{\partial}{\partial Z}
- \epsilon \sum_{k=1}^{2}
C_{i3kl}(x,Z) \frac{\partial}{\partial Z}
\frac{\partial}{\partial x_k}
- \epsilon \sum_{k=1}^{2}
\left( \frac{\partial}{\partial Z} C_{i3kl}(x,Z) \right)
\frac{\partial}{\partial x_k}\\
- \epsilon^2 \sum_{j,k=1}^{2}
C_{ijkl}(x,Z)\frac{\partial}{\partial x_j}
\frac{\partial}{\partial x_k}-\epsilon \sum_{j=1}^{2}
\left(\frac{\partial}{\partial x_j}C_{ij3l}(x,Z) \right)
\frac{\partial}{\partial Z}\\
- \epsilon^2 \sum_{j,k=1}^{2}
\left(\frac{\partial}{\partial x_j}C_{ijkl}(x,Z) \right)
\frac{\partial}{\partial x_k}.
\end{multline}
With Definition~\ref{weyl} in Appendix~A, we consider $\hat{H}$ as a semiclassical pseudodifferential operator on $\mathbb{R}^2$, where $x$ belongs to; then $(\ref{semi1})$ can be rewritten as
\begin{equation}\label{sharp}
\epsilon^2\partial_t^2 v+\text{Op}_\epsilon(H(x,\xi))v\ \sim\ 0,
\end{equation}
where $\text{Op}_\epsilon(H(x,\xi))$ is a semiclassical pseudodifferential operator with symbol $H(x,\xi)$ defined by
\[H(x,\xi)=H_0(x,\xi)+H_1(x,\xi)\]
with
\begin{multline}\label{H0}
H_{0,il}(x,\xi) =
-\frac{\partial}{\partial Z}C_{i33l}(x,Z)
\frac{\partial}{\partial Z}
\\
- \ii\sum_{j=1}^{2}
C_{ij3l}(x,Z) \xi_j
\frac{\partial}{\partial Z}
- \ii\sum_{k=1}^{2}
C_{i3kl}(x,Z) \frac{\partial }{\partial Z}
\xi_k
- \ii \sum_{k=1}^{2}
\left( \frac{\partial}{\partial Z} C_{i3kl}(x,Z) \right)
\xi_k 
\\
+ \sum_{j,k=1}^{2}
C_{ijkl}(x,Z) \xi_j \xi_k
\hspace*{4.0cm}
\end{multline}
and
\begin{equation}\label{H1}
H_{1,il}(x,\xi) =-\epsilon \sum_{j=1}^{2}
\left(\frac{\partial}{\partial x_j}C_{ij3l}(x,Z) \right)
\frac{\partial}{\partial Z}
- \ii\epsilon \sum_{j,k=1}^{2}
\left(\frac{\partial}{\partial x_j}C_{ijkl}(x,Z) \right)
\xi_k.
\end{equation}
We view $H_0(x,\xi)$ and $H_1(x,\xi)$ as ordinary differential
operators in $Z$, with domain
\[
   \mathcal{D} = \left\{ v \in H^2(\mathbb{R}^-)\ \bigg|\
      \sum_{l=1}^3\left(C_{i33l}(x,0)
    \frac{\partial v_l}{\partial Z}(0)
       + \ii \sum_{k=1}^2C_{i3kl}\xi_kv_l(0)\right) = 0 \right\} .
\]

\subsection{Effective Hamiltonians}

We use eigenvalues and eigenfunctions of $H_0(x,\xi)$ to construct
approximate solutions to $(\ref{sharp})$. For fixed
$(x,\xi)$, an eigenvalue $\Lambda(x,\xi)$ and the corresponding
eigenfunction $V(x,\xi)$ of $H_0(x,\xi)$ are such that
\begin{equation}\label{eigeq}
   H_0(x,\xi)V(x,\xi) = V(x,\xi) \Lambda(x,\xi),
\end{equation}
where $V(x,\xi)\in \mathcal{D}$. Since $H_0(x,\xi)$ is a positive
symmetric operator in $L^2(\mathbb{R}^-)$ for $\xi\neq 0$,
$\Lambda(x,\xi)$ is real-valued, and also positive. We note, here,
that the contribution coming from $\xi$ confined to a compact set is
negligible. This is why we can assume that $\xi \ne 0$.

\medskip

We let $\mathcal{L}(\mathcal{H}_1,\mathcal{H}_2)$ denote the set of
all bounded operators from a normed space $\mathcal{H}_1$ to a normed
space $\mathcal{H}_2$.

\medskip\medskip

\begin{theorem}\label{conjugation_sys}
Let $\Lambda_\alpha(x,\xi)$ be an eigenvalue of $H_0(x,\xi)$, and $U
\subset T^* \mathbb{R}^2\setminus 0$ be open. Assume that
$\Lambda_\alpha(x,\xi)$ has constant multiplicity $m_\alpha$ for all
$(x,\xi)\in U$. There exist
$\Phi_{\alpha,m}(x,\xi)\in\mathcal{L}(\mathcal{D},L^2(\mathbb{R}^-))$
and $a_{\alpha,m}(x,\xi)\in\mathcal{L}(L^2(\mathbb{R}^-),\mathcal{D})$
which admits asymptotic expansions
\begin{equation}\label{eq:symb_sys}
\begin{split}
\Phi_{\alpha,\epsilon}(x,\xi) \sim
\sum_{m=0}^{\infty} \Phi_{\alpha,m}(x,\xi) \, \epsilon^m ,
\\
a_{\alpha,\epsilon}(x,\xi) \, \sim \, \sum_{m=0}^{\infty} a_{\alpha,m}(x,\xi) \, \epsilon^m ,
\end{split}
\end{equation}
and satisfy
\begin{equation}\label{eq:comm_sys}
    H \circ \Phi_{\alpha,\epsilon}(x,\xi)\, 
   = \Phi_{\alpha,\epsilon}
    \circ a_{\alpha,\epsilon}(x,\xi)+\mathcal{O}(\epsilon^\infty)\,
\end{equation}
where $\circ$ denotes the composition of symbols (see Appendix~A).
Furthermore, $a_{\alpha,0}(x,\xi)=\Lambda_\alpha(x,\xi)I$ and
$\Phi_{\alpha,0}(x,\xi)$ is the projection onto the eigenspace
associated with $\Lambda_\alpha(x,\xi)$.
\end{theorem}

\medskip\medskip

\begin{proof}
First, we note that $\Lambda_\alpha(x,\xi)\in C^\infty(U)$ (cf. Appendix B of \cite{CDV3}). By the composition of symbols (cf. Appendix~A), we have
\[H\circ\Phi_{\alpha,\epsilon}=H_0\Phi_{\alpha,0}+\epsilon H_0\Phi_{\alpha,1}+\epsilon H_1\Phi_{\alpha,0}+\frac{\epsilon}{2\ii}\sum_{j=1}^2\frac{\partial H_0}{\partial \xi_j}\frac{\partial \Phi_{\alpha,0}}{\partial x_j}+\mathcal{O}(\epsilon^2)\]
and 
\[\Phi_{\alpha,\epsilon}\circ a_{\alpha,\epsilon} =\Phi_{\alpha,0}\Lambda_\alpha+\epsilon\Phi_{\alpha,0}a_{\alpha,1}+\epsilon\Phi_{\alpha,1}\Lambda_\alpha+\frac{\epsilon}{2\ii}\sum_{j=1}^2\frac{\partial \Phi_{\alpha,0}}{\partial \xi_j}\frac{\partial \Lambda_{\alpha,0}}{\partial x_j}+\mathcal{O}(\epsilon^2).\]
We construct the terms $\Phi_{\alpha,m}$ by collecting terms of equal orders in the two expansions above. Terms of order $\mathcal{O}(\epsilon^0)$ give
\[H_0\Phi_{\alpha,0}=\Phi_{\alpha,0}\Lambda_\alpha,\]
which is consistent with (\ref{eigeq}). Terms of order $\mathcal{O}(\epsilon^1)$ give
\[(H_0(x,\xi)-\Lambda_\alpha(x,\xi))\Phi_{\alpha,1}(x,\xi)=\Phi_{\alpha,0}a_{\alpha,1}+R(\Phi_{\alpha,0}),\]
where $R(\Phi_{\alpha,0})\in\mathcal{L}(\mathcal{D},L^2(\mathbb{R}^-))$ denotes the remaining terms.
We choose $a_{\alpha,1} $, so that, for  $u\in\mathcal{D}$,
$(\Phi_{\alpha,0}a_{\alpha,1}+R(\Phi_{\alpha,0}))u$ is orthogonal
to the eigenspace of $\Lambda_\alpha(x,\xi)$. Then we let
$v=\Phi_{\alpha,1}(x,\xi)u$ be the unique solution of
\[(H_0(x,\xi)-\Lambda_\alpha(x,\xi))v=(\Phi_{\alpha,0}a_{\alpha,1}+R(\Phi_{\alpha,0}))u ,\]
which is orthogonal to the eigenspace of
$\Lambda_\alpha(x,\xi)$. Thus we have defined the operator
$\Phi_{\alpha,1}(x,\xi)$. Higher order terms $\Phi_{\alpha,m}(x,\xi)$,
$m \ge 2$, can be constructed successively by solving the equations,
\[
   (H_0(x,\xi)-\Lambda_\alpha(x,\xi)) \Phi_{\alpha,m}(x,\xi)
     = R(\Phi_{\alpha,0},\cdots,\Phi_{\alpha,m-1}),
\]
with $R(\Phi_{\alpha,0},\cdots,\Phi_{\alpha,m-1})
\in\mathcal{L}(\mathcal{D}, L^2(\mathbb{R}^-))$, since
$\Phi_{\alpha,k},a_{\alpha,k}$, $0\leq k\leq m-1$, in their respective
topologies, depend continuously on $(x,\xi)$ by induction. Hence, we
can choose \\ $a_{\alpha,m} =
-R(\Phi_{\alpha,0},\cdots,\Phi_{\alpha,m-1})$ and solve for
$\Phi_{\alpha,m}$. This completes the construction of $a_{\alpha,m}$
and $\Phi_{\alpha,m}$.
\end{proof}

\medskip\medskip

For any bounded set $U' \subset T^*\mathbb{R}^2 \setminus 0$, we let
$U \subset T^*\mathbb{R}^2$ be such that $\overline{U'} \subset U$ and
$\chi \in C_0^\infty(U)$ with $\chi \equiv 1$ on $\overline{U'}$. Then
\begin{equation}\label{defchi}
   J_{\alpha,\epsilon} = \text{Op}_\epsilon
      \left(\frac{1}{\sqrt{\epsilon}} \chi(x,\xi)
                         \Phi_{\alpha,\epsilon}\right)
\end{equation}
defines a linear map on $\mathcal{M}(V)$ (the factor $\frac{1}{\sqrt{\epsilon}}$ is to make the operator $J_{\alpha,\epsilon}$ (microlocally) unitary), where
\[
   V = \{(x,\xi,Z,\zeta) \in T^* X\ :\ (x,\xi) \in U' \}
\]
and $\mathcal{M}(V)$ denotes the space of microfunctions on $V$; see,
Definition \ref{A4} of Appendix~\ref{AA}. Moreover, we find that, using
Theorem~\ref{conjugation_sys},
\begin{equation}
   \hat{H} J_{\alpha,\epsilon}
        = J_{\alpha,\epsilon} \text{Op}_\epsilon(a_{\alpha,\epsilon}),
\end{equation}
where the left- and right-hand sides are considered as microfunctions.
~\\

\begin{remark}
We note that the principal symbol of $a_{\alpha,\epsilon}$ is
real-valued. Hence, the propagation of wavefront set of surface waves
is purely on the the surface. However, since the lower order terms of
$a_{\alpha,\epsilon}$ are operators in $Z$, there exists coupling of
surface waves with interior motion. There is no such coupling if the
eigenvalue $\Lambda_\alpha$ has multiplicity one. For this case
$a_{\alpha,\epsilon}$ is a classical pseudodifferential symbol.
\end{remark}

\subsection{Surface-wave equations and parametrices}\label{param}

With the results of the previous subsection, we construct approximate
solutions of the system \eqref{elaswaeq} with initial values
\[
h(x,\epsilon Z)=\sum_{\alpha=1}^{\mathfrak{M}}
	J_{\alpha,\epsilon} W_{\alpha}(x,Z).
\]
Then we construct solutions of the form
\begin{equation}\label{ans-wa_sys}
\begin{split}
   u(t,x,z) = v\left(t,x,\frac{z}{\epsilon}\right)
   &= \sum_{\alpha=1}^{\mathfrak{M}}
   v_{\alpha,\epsilon}\left(t,x,\frac{z}{\epsilon}\right) ,
\\
   v_{\alpha,\epsilon}\left(t,x,\frac{z}{\epsilon}\right)
   &= J_{\alpha,\epsilon} W_{\alpha,\epsilon}(t,x,\frac{z}{\epsilon}) ,
\end{split}
\end{equation}
with
\begin{eqnarray*}
	v\left(0,x,\frac{z}{\epsilon}\right) = 0 ,
	\quad
	\p_t v\left(0,x,\frac{z}{\epsilon}\right)
	&=& \sum_{\alpha=1}^{\mathfrak{M}}
	J_{\alpha,\epsilon} W_{\alpha}(x,\frac{z}{\epsilon}) ,
\\
	\frac{c_{i3kl}}{\rho}
        \partial_k v_l\left(t,x,\frac{z}{\epsilon}\right)
                   \bigg|_{z=0} &=& 0.
\end{eqnarray*}
Here, $\mathfrak{M}$ is chosen such that for each $(x,\xi)\in U$, there
are at least $\mathfrak{M}$ eigenvalues for $H_0(x,\xi)$. We assume that
all eigenvalues $\Lambda_1 < \cdots <
\Lambda_\alpha<\cdots<\Lambda_{\mathfrak{M}}$ are of constant
multiplicities, $m_1,\cdots,m_\alpha,\cdots,m_{\mathfrak{M}}$. We let
$W_{\alpha,\epsilon} $ solve the initial value problems
\begin{eqnarray}
   [\epsilon^2 \p_t^2
       + \operatorname{Op}_{\epsilon}(a_{\alpha,\epsilon})(.,D_x)]
	 W_{\alpha,\epsilon}(t,x,Z) &=& 0 ,
\label{wqeaj_sys}\\
   W_{\alpha,\epsilon}(0,x,Z) &=& 0 ,\quad
   \p_t W_{\alpha,\epsilon}(0,x,Z) = J_{\alpha,\epsilon}W_{\alpha}(x,Z) ,
\label{eq:initrange}
\end{eqnarray}
$\alpha = 1,\ldots,\mathfrak{M}$. Equation (\ref{eq:initrange}) means
that the initial values are in the span of the ranges of operators
$\op_\epsilon(\Phi_{\alpha,\epsilon})$, $\alpha=
1,\ldots,\mathfrak{M}$.

\medskip

\medskip

We address the existence of eigenvalues of $H_0(x,\xi)$ in later
sections, under the Assumption \ref{assu_tensor}. To construct the
parametrix, we use the first-order system for $W_{\epsilon}$ that is
equivalent to (\ref{wqeaj_sys}),
\begin{equation} \label{eq:first_order_time}
	\epsilon\pdpd{}{t} \left( \begin{array}{cc} W_{\alpha,\epsilon} \\ \epsilon\pdpd{W_{\alpha,\epsilon}}{t}
	\end{array} \right)
	= \left( \begin{array}{cc} 0 & 1 \\ -\operatorname{Op}_{\epsilon}(a_{\alpha,\epsilon})
		& 0 \end{array} \right)
	\left( \begin{array}{cc} W_{\alpha,\epsilon} \\ \epsilon\pdpd{W_{\alpha,\epsilon}}{t} \end{array} \right) .
\end{equation}
We write
$\operatorname{Op}_{\epsilon}(b_{\alpha,\epsilon}) =
[\operatorname{Op}_{\epsilon}(a_{\alpha,\epsilon})]^{1/2}$. The principal
symbol of $\operatorname{Op}_{\epsilon}(b_{\alpha,\epsilon})$ is given by
$\Lambda^{1/2}_\alpha(x,\xi)$. Then
\begin{equation}
	\label{eq:decoupled_firstorder_u_f}
	W_{\alpha,\epsilon,\pm} = \half W_{\alpha,\epsilon} \pm
	\half \ii
	\epsilon\operatorname{Op}_{\alpha,\epsilon}(b_{\alpha,\epsilon})^{-1} \pdpd{W_{\alpha,\epsilon}}{t} ,
\end{equation}
satisfy the two first-order (``half wave'') equations
\begin{equation} \label{eq:decoupled_firstorder}
	P_{\alpha,\epsilon,\pm}(x,D_x,D_t) W_{\alpha,\epsilon,\pm} = 0 ,
\end{equation}
where
\begin{equation}
	P_{\alpha,\epsilon,\pm}(x,D_x,D_t)
	= \epsilon\p_t \pm \ii \operatorname{Op}_{\epsilon}(b_{\alpha,\epsilon}) ,
\end{equation}
supplemented with the initial conditions
\begin{equation}\label{eq:hpolIVs}
   h_{\alpha,\epsilon,\pm} = \pm \half
   \ii \operatorname{Op}_{\alpha,\epsilon}(b_{\alpha,\epsilon})^{-1}J_{\alpha,\epsilon}
            W_{\alpha,\epsilon} |_{t=0} .
\end{equation}
The principal symbol of $\frac{1}{\ii}P_{\alpha,\epsilon,\pm}=\epsilon D_t
\pm\text{Op}_\epsilon(b_{\alpha,\epsilon})$ is given by $\omega \pm
\Lambda^{1/2}_\alpha(x,\xi)$ defining a Hamiltonian. The solution to the
original problem is then given by
\begin{equation}\label{summ}
   W_{\alpha,\epsilon} = W_{\alpha,\epsilon,+} + W_{\alpha,\epsilon,-} .
\end{equation}

We construct WKB solutions for operator $P_{\epsilon,+}(x,D_x,D_t)$,
introducing a representation of the solution operator of the form,
\begin{equation}\label{FIO}
   S_{\alpha,\epsilon}(t) h_{\alpha,\epsilon,+}(x) := \frac{1}{(2\pi\epsilon)^2}
     \iint e^{\frac{\ii (\psi_{\alpha,+}(x,\eta,t) - \langle y,\eta \rangle)}{
       \epsilon}} \mathcal{A}_{\alpha,+}(x,\eta,t,\epsilon) \,
     h_{\alpha,\epsilon,+}(y) \mathrm{d}y \mathrm{d}\eta.
\end{equation}
For more details of this construction, we refer to
\cite{Bardos,Zworski}.

\medskip\medskip

\noindent
\textbf{Construction of the phase function}. Surface-wave
bicharacteristics are solutions, $(x,\xi)$, to the Hamiltonian system,
\begin{equation}\label{bichar}
   \frac{\partial x_k(y,\eta,t)}{\partial t}
       = \frac{\partial\Lambda^{1/2}_\alpha(x,\xi)}{
                \partial \xi_k} ,\
   \frac{\partial \xi_k(y,\eta,t)}{\partial t}
       =-\frac{\partial\Lambda^{1/2}_\alpha(x,\xi)}{
                \partial x_k} ,\quad
    (x,\xi)|_{t=0} = (y,\eta) .
\end{equation}
We define
\[
   \kappa_t(y,\eta)=(x,\xi) .
\]
For $(y,\eta)$ in $U$, the map
\[
   (x,\xi,y,\eta) = (\kappa_t(y,\eta),(y,\eta)) \mapsto (x,\eta)
\]
is surjective for small $t$. We write $\xi=\Xi(x,\eta,s)$ and
introduce the phase function
\begin{equation}
   \psi_{\alpha,+}(x,\eta,t) = \langle
x, \eta\rangle-\int_0^t\Lambda^{1/2}_\alpha(x,\Xi(x,\eta,s))\mathrm{d}s,
\end{equation}
where $\langle x,\eta\rangle$ denotes the inner product of $x$ and $\eta$.
This phase function solves
\[
   \partial_t \psi_{\alpha,+} + \Lambda^{1/2}_\alpha(x,\partial_x\psi_{\alpha,+}) = 0 ,~~
   \psi_{\alpha,+}(x,\eta,0)=\langle x, \eta\rangle
\]
\cite[Chapter~10]{Zworski}. The group velocity, $\mathfrak{v}$, is
defined by
\begin{equation}\label{groupvelocity}\mathfrak{v}_k:=\frac{\partial x_k}{\partial t}=\frac{\partial\Lambda^{1/2}_\alpha(x,\xi)}{\partial \xi_k}.\end{equation}
We note that $\Lambda^{1/2}_\alpha(x,\xi)$ is non-homogeneous in $\xi$,
thus the group velocity and the bicharacteristics depend on $|\xi|$.

\medskip\medskip

\noindent\textbf{Construction of the amplitude}. The amplitude
$\mathcal{A}$ in (\ref{FIO}) must satisfy 
\[
   (\epsilon D_t + b_{\alpha,\epsilon}(x,\epsilon D))
        (e^{\frac{\ii\psi_{\alpha,+}}{\epsilon}} \mathcal{A}_{\alpha,+})
                        = \mathcal{O}(\epsilon^\infty) ,
\]
or, equivalently,
\begin{equation}\label{asym1}
   (\partial_t \psi_{\alpha,+} + \epsilon D_t
        + e^{-\frac{\ii\psi_{\alpha,+}}{\epsilon}} b_{\alpha,\epsilon}(x,\epsilon D)
     e^{\frac{\ii\psi_{\alpha,+}}{\epsilon}})\mathcal{A}_{\alpha,+}
            = \mathcal{O}(\epsilon^\infty) .
\end{equation}
We note that $\eta$ appears in the construction of the phase function
$\psi_{\alpha,+}$. We construct solutions, $\mathcal{A}_{\alpha,+}$, via the expansion,
\[
   \mathcal{A}_{\alpha,+}(x,\eta,t) \sim
       \sum_{k=0}^\infty \epsilon^k \mathcal{A}_{\alpha,+,k}(x,\eta,t) ,
\]
and writing the symbol of $b_{\alpha,\epsilon}(x,\epsilon D)$ as an expansion,
\[
   b_{\alpha,\epsilon}(x,\xi) = \sum_{j=0}^\infty\epsilon^j b_{\alpha,j}(x,\xi) .
\]
By the asymptotics of
$e^{-\frac{\ii\psi_{\alpha,+}}{\epsilon}}b_{\alpha,\epsilon}(x,\epsilon D)
e^{\frac{\ii\psi_{\alpha,+}}{\epsilon}}$, the terms of order
$\mathcal{O}(\epsilon^0)$ yield the equation,
\begin{equation}\label{phase_eq}
   \partial_t \psi_{\alpha,+} + b_{\alpha,0}(x,\partial_x\psi_{\alpha,+}) = \partial_t \psi_{\alpha,+}
           + \Lambda^{1/2}_\alpha(x,\partial_x\psi_{\alpha,+}) = 0 ,
\end{equation}
for the phase function $\psi_{\alpha,+}$, which was introduced above. The terms
of order $\mathcal{O}(\epsilon^1)$ yield the transport equation for
$\mathcal{A}_{\alpha,+,0}$, 
\[
\begin{split}
   (D_t + L) \mathcal{A}_{\alpha,+,0} &= 0 ,\\[0.25cm]
   \mathcal{A}_{\alpha,+,0}(x,\eta,0) &= 1 ,
\end{split}
\]
where
\begin{equation}\label{L}
   L = \sum_{j=1}^2 \partial_{\xi_j}
        \Lambda^{1/2}_\alpha(x,\partial_x\psi_{\alpha,+})D_{x_j}
       + b_{\alpha,1}(x,\partial_x\psi_{\alpha,+}) - \frac{\ii}{2} \sum_{j,k=1}^2
    \partial_{\xi_k\xi_j}\Lambda^{1/2}_\alpha(x,\partial_x\psi_{\alpha,+})D_{x_kx_l}\psi_{\alpha,+} ,
\end{equation}
in which $\partial_{\xi_k\xi_j} = \partial_{\xi_k}\partial_{\xi_j}$,
$D_{x_kx_l} = D_{x_k}D_{x_l}$. Denote
\[c_{t,\epsilon}(x,\epsilon D,\eta)=e^{-\frac{\ii\psi_{\alpha,+}(x,\eta,t)}{\epsilon}}b_{\alpha,\epsilon}(x,\epsilon D)
e^{\frac{\ii\psi_{\alpha,+}(x,\eta,t)}{\epsilon}}=\sum_{l=0}^\infty\epsilon^lc_{t,l}(x,\epsilon D,\eta).\]

We can construct $\mathcal{A}_{\alpha,+,k}$, $k =
1,2,\cdots$, recursively by solving transport equations of the form
\cite{Bardos}
\[
\begin{split}
   (D_t + L) \mathcal{A}_{\alpha,+,k} &= F(\mathcal{A}_{\alpha,+,0},\cdots,\mathcal{A}_{\alpha,+,k-1}) ,
   ~~k \geq 1 ,\\[0.25cm]
   \mathcal{A}_{\alpha,+,k}(x,\eta,0) &= 0 .
\end{split}
\]
Here,
\[F(\mathcal{A}_{\alpha,+,0},\cdots,\mathcal{A}_{\alpha,+,k-1})=-\ii\partial\psi_{\alpha,+}\mathcal{A}_{\alpha,+,k-1}(x,\eta,t)-\ii \sum_{j+l=k-1}c_{t,l}(x,\epsilon D,\eta)\mathcal{A}_{\alpha,+,j}(x,\eta,t).\]

For solving the transport equations, we use the bicharacteristics
determined by (\ref{bichar}). We note that
\[
   D_t + \sum_{j=1}^2 \partial_{\xi_j}
        \Lambda^{1/2}_\alpha(x,\partial\psi_{\alpha,+})D_{x_j} = D_t
   + \sum_{j=1}^2 \frac{\partial x_j(y,\eta,t)}{\partial t} D_{x_j}
\]
so that
\[
\begin{split}
\left(\frac{\mathrm{d}}{\mathrm{d}t}+\ii b_{\alpha,1}(x,\partial_x\psi_{\alpha,+})+\frac{1}{2}\sum_{j,k=1}^2\partial_{\xi_k\xi_j}\Lambda^{1/2}_\alpha(x,\partial_x\psi_{\alpha,+})D_{x_kx_l}\psi_{\alpha,+}\right) \mathcal{A}_{\alpha,+,k}(x(y,\eta,t),\eta,t)\\
~~~~~~~~~~~~~=F(\mathcal{A}_{\alpha,+,0},\cdots,\mathcal{A}_{\alpha,+,k-1})\\
\mathcal{A}_{\alpha,+,0}(x(y,\eta,0),\eta,0)=1,~~~~~~~\mathcal{A}_{\alpha,+,k}(x(y,\eta,0),\eta,0)=0,~~~\text{for}~k= 1,2,3,\cdots.
\end{split}\]
The equations above can be solved using the standard theory of
ordinary differential equations. For more details about semiclassical Fourier Integral Operators, we refer to Chapter 10 of Zworski's book \cite{Zworski}.

\medskip\medskip

\noindent
\textbf{Trace}. We can write the solution to
(\ref{eq:decoupled_firstorder}), up to leading order, as
\[
\frac{1}{(2\pi\epsilon)^2}\iint \mathcal{G}_{\alpha,\pm,0}(x,t,\epsilon)h_{\epsilon,\pm}(y) \exp(-\frac{\mathrm{i}}{\epsilon}\langle x-y,\eta\rangle)\mathrm{d}y\mathrm{d}\eta
\]
with
\[
\mathcal{G}_{\alpha,\pm,0}(x,t,\eta,\epsilon)=\exp\left[ \mp\frac{\mathrm{i}}{\epsilon}\int_0^t\Lambda^{1/2}_\alpha(x,\Xi(x,\eta,s))\mathrm{d}s\right]\mathcal{A}_{\alpha,\pm,0}(x,\eta,t).
\]
Then, using (\ref{eq:decoupled_firstorder_u_f}), (\ref{eq:hpolIVs})
and (\ref{summ}), we can write down the approximate Green's function
(microlocalized in $x$), still up to leading order, for surface waves
\begin{multline}
\mathcal{G}_0(Z,x,t,Z',\eta;\epsilon)=\sum_{\alpha=1}^{\mathfrak{M}} J^{(x)}_{\alpha,\epsilon}(Z,x,\eta)
\\
\left(\frac{\mathrm{i}}{2}\mathcal{G}_{\alpha,+,0}(x,t,\eta,\epsilon)-\frac{\mathrm{i}}{2}\mathcal{G}_{\alpha,-,0}(x,t,\eta,\epsilon)\right)\Lambda^{-1/2}_\alpha(x,\eta)J^{(x)}_{\alpha,\epsilon}(Z',x,\eta)
\\
=\sum_{\alpha=1}^{\mathfrak{M}} J^{(x)}_{\alpha,\epsilon}(Z,x,\eta)
\Bigg(\frac{\mathrm{i}}{2}\exp\left[-\frac{\mathrm{i}}{\epsilon}\int_0^t\Lambda^{1/2}_\alpha(x,\Xi(x,\eta,s))\mathrm{d}s\right]\mathcal{A}_{\alpha,+,0}(x,\eta,t)
\\
-\frac{\mathrm{i}}{2}\exp\left[\frac{\mathrm{i}}{\epsilon}\int_0^t\Lambda^{1/2}_\alpha(x,\Xi(x,\eta,s))\mathrm{d}s\right]\mathcal{A}_{\alpha,-,0}(x,\eta,t)\Bigg)\Lambda^{-1/2}_\alpha(x,\eta)J^{(x)}_{\alpha,\epsilon}(Z',x,\eta).
\end{multline}

Then
\begin{equation*}
\begin{split}
&\epsilon\partial_t\mathcal{G}_0(Z,x,t,Z',\eta;\epsilon)\\
=&\sum_{\alpha=1}^{\mathfrak{M}} J^{(x)}_{\alpha,\epsilon}(Z,x,\eta)
\Bigg(\frac{1}{2}\exp\left[-\frac{\mathrm{i}}{\epsilon}\int_0^t\Lambda^{1/2}_\alpha(x,\Xi(x,\eta,s))\mathrm{d}s\right]\mathcal{A}_{\alpha,+,0}(x,\eta,t)
\\
&+\frac{1}{2}\exp\left[\frac{\mathrm{i}}{\epsilon}\int_0^t\Lambda^{1/2}_\alpha(x,\Xi(x,\eta,s))\mathrm{d}s\right]\mathcal{A}_{\alpha,-,0}(x,\eta,t)\Bigg)\Lambda^{1/2}_\alpha(x,\Xi(x,\eta,t))\Lambda^{-1/2}_\alpha(x,\eta)J^{(x)}_{\alpha,\epsilon}(Z',x,\eta).
\end{split}
\end{equation*}

Now, we take the semiclassical Fourier transform in $t,$ multiply by the test function $\chi(\omega)\in \mathcal{S}(\R)$ and apply the stationary phase formula to the leading order (see, for example, \cite{Dimassi-Sjostrand}, Proposition 5.2). We get
\begin{align*}&\int\widehat{\epsilon\partial_t\mathcal{G}_0}(Z,x,\omega,Z',\eta;\epsilon)\chi(\omega)d\omega=\frac{1}{2\pi\epsilon}\int\!\!\int \epsilon\partial_t\mathcal{G}_0(Z,x,t,Z',\eta;\epsilon)e^{-\frac{\mathrm{i}t\omega}{\epsilon}}\mathrm{d}t\chi(\omega)d\omega\sim\sum_{\alpha=1}^{\mathfrak{M}} J^{(x)}_{\alpha,\epsilon}(Z,x,\eta)\\
&
\Bigg(\frac{1}{2}\mathcal{A}_{\alpha,+,0}(x,\eta,0)\chi(-\Lambda^{1/2}_\alpha(x,\eta))+\frac{1}{2}\mathcal{A}_{\alpha,-,0}(x,\eta,0)\chi(\Lambda^{1/2}_\alpha(x,\eta))\Bigg)J^{(x)}_{\alpha,\epsilon}(Z',x,\eta).
\end{align*}

As $\mathcal{A}_{\alpha,\pm,0}(x,\eta,0)=1$, we get
\begin{align*}&\widehat{\epsilon\partial_t\mathcal{G}_0}(Z,x,\omega,Z',\eta;\epsilon)\\ \sim&\sum_{\alpha=1}^{\mathfrak{M}} \frac{1}{2}J^{(x)}_{\alpha,\epsilon}(Z,x,\eta)\left[\delta(\omega+\Lambda_\alpha^{1/2}(x,\eta))+\delta(\omega-\Lambda_\alpha^{1/2}(x,\eta))\right]J^{(x)}_{\alpha,\epsilon}(Z',x,\eta)\\
=&\sum_{\alpha=1}^{\mathfrak{M}} J^{(x)}_{\alpha,\epsilon}(Z,x,\eta)\delta\left(\omega^2-\Lambda_\alpha(x,\eta)\right)\Lambda^{1/2}_\alpha(x,\eta)J^{(x)}_{\alpha,\epsilon}(Z',x,\eta)
\end{align*}
On the diagonal $Z=Z'$,  we get
$$ \widehat{\epsilon\partial_t\mathcal{G}_0}(Z,x,\omega,Z,\eta;\epsilon)\sim \sum_{\alpha=1}^{\mathfrak{M}} \left(J^{(x)}_{\alpha,\epsilon}(Z,x,\eta)\right)^2\delta\left(\omega^2-\Lambda_\alpha(x,\eta)\right)\Lambda^{1/2}_\alpha(x,\eta).$$

We remember here that $J^{(x)}_{\alpha,\epsilon}(Z,x,\eta)=
\frac{1}{\sqrt{\epsilon}}(\Phi_{\alpha,0}(Z,x,\eta)+{\mathcal O}(\epsilon^{-1}))$ for $(x,\eta)\subset U'$ (cf. (\ref{defchi})). Assuming the eigenvalues are all simple, we can simply take $\Phi_{\alpha,0}$ as normalized eigenfunctions. Then, 
\[
\int_{\mathbb{R}^-}\widehat{\epsilon\partial_t\mathcal{G}_0}(Z,x,\omega,Z,\eta;\epsilon)\mathrm{d}\epsilon Z= \sum_{\alpha=1}^{\mathfrak{M}}\delta(\omega^2-\Lambda_\alpha(x,\eta))\Lambda^{1/2}_\alpha(x,\eta)+{\mathcal O}(\epsilon^{-1}).
\]

The left-hand side can be viewed as a local trace in the boundary
normal coordinate.

\medskip\medskip

\begin{remark}\label{rem}
In the oscillatory integral $(\ref{FIO})$, the phase variables are the
components of $\eta$. We can construct an alternative representation by
using the frequency, $\omega$, as a phase variable. This would result in the following form of solution operator
\begin{equation}\label{FIOw}
   S_{\alpha,\epsilon}(t) h_{\alpha,\epsilon,+}(x) := \frac{1}{(2\pi\epsilon)^2}
     \iint e^{\frac{\ii (\phi_\alpha(x,x_0,\omega) - \langle\omega,t \rangle)}{
       \epsilon}} \mathcal{A}_{\alpha,+}(x,x_0,\omega,\epsilon) \,
     h_{\alpha,\epsilon,+}(x_0) \mathrm{d}x_0 \mathrm{d}\omega.
\end{equation}
In a neighborhood of
a fixed point, $x_0$, we construct $\phi_\alpha = \phi_\alpha(x,x_0,\omega)$ satisfying
the eikonal equation
\begin{equation}\label{eikonal}
   \Lambda^{1/2}_\alpha(x,\partial_x\phi_\alpha(x,x_0,\omega))\ =\ \omega .
\end{equation}
Then the phase function is given by
\[
   \tilde{\psi}_\alpha(x,x_0,t,\omega) = \phi_\alpha(x,x_0,\omega) - \omega t
\]
and satisfies $(\ref{phase_eq})$,
\[
   \partial_t \tilde{\psi}_\alpha
         + \Lambda^{1/2}_\alpha(x,\partial_x \tilde{\psi}_\alpha) = 0 .
\]
We consider the bicharacteristics $(y(x_0,\xi_0,t),\eta(x_0,\xi_0,t))$
and let $(x,\xi) = (y,\eta)|_{t=t_x}$, while $\partial_x
\phi_\alpha(x,x_0,\omega)$ \\ $= \xi$. The generating function $\phi_\alpha$ satisfies
\[
   \frac{\partial}{\partial t}\phi_\alpha(y(t),x_0,\omega)
            = \langle\eta(x_0,\xi_0,t)\, ,\,\dot{y}(x_0,\xi_0,t)\rangle,
\]
and can then be written in the form 
\[
   \phi_\alpha(x,x_0,\omega) = \int_{0}^{t_x} \langle\eta(x_0,\xi_0,t)
         \, ,\, \dot{y}(x_0,\xi_0,t)\rangle\mathrm{d}t
   = \int_{x_0}^x \sum_{j=1}^2 \eta_j \mathrm{d}y_j .
\]

We mention an alternative representation. For any $\hat{\xi}_0 \in
S^2$, we let $\xi_0 = K_0 \hat{\xi}_0$, $K_0(x_0,\omega,\hat{\xi}_0)$
be the solution of
\[
   \Lambda^{1/2}_\alpha(x_0,K_0\hat{\xi}_0) = \omega .
\]
This $K_0$ is unique due to the monotonicity of eigenvalue in $|\xi|$
with respect to the quadratic form associated with $H_0$. We denote
$\eta(x_0,\xi_0,t) = K(x_0,\xi_0,t) \hat{\eta}$ with $\hat{\eta} =
\frac{\eta}{|\eta|} \in S^2$, and define the phase velocity,
$\mathfrak{V} = \mathfrak{V}(y,\omega,\hat{\eta})$, as
\begin{equation}\label{phasevelocity}
   \mathfrak{V} = \frac{\omega}{K} .
\end{equation}
Thus $K(x_0,\xi_0,0)=K_0$. Using $(\ref{groupvelocity})$ and
$(\ref{phasevelocity})$ we find that
\begin{equation}
   \phi_\alpha(x,x_0,\omega)
       = \omega \, \int_0^{t_x} \sum_{j=1}^2 \hat{\eta}_j
         \frac{\mathfrak{v}_j}{\mathfrak{V}} \,  \mathrm{d}t ,
\end{equation}
defined in a neighborhood of $x_0$. Since
$\frac{\mathfrak{v}}{\mathfrak{V}}$ is frequency dependent, the
geodesic distance is frequency dependent.
\end{remark}

\section{Characterization of the spectrum of $H_0$}
\label{sec:3}

In this section, we characterize the spectrum of $H_0$. In the case of
constant coefficients, that is, a homogeneous medium, the spectrum was
described in \cite{Kamo}. In the case of isotropic media, the spectrum
was studied by Colin de Verdi\`ere \cite{Yves1}.

We view $$h_0(x,\hat{\xi}) = \frac{1}{|\xi|^2} H_0(x,\xi)$$ as a
semiclassical pseudodifferential operator with $\frac{1}{|\xi|}$ as
the semiclassical parameter, $h$ say. We write
$\hat{\xi}=\frac{\xi}{|\xi|}$ and use $\zeta$ to denote the Fourier
variable for $\frac{1}{\ii |\xi|}\frac{\partial}{\partial Z}$. Then
the principal symbol of $h_0(x,\hat{\xi})$ is given by
\begin{equation}\label{principle}
   h_0(x,\hat{\xi})(Z,\zeta) = T \zeta^2 + (R + R^T) \zeta + Q ,
\end{equation}
where
\begin{equation}\label{TRQ}
\begin{split}
   & T _{il}(x,Z):= C_{i33l}(x,Z) ,\\
    & R_{il}(x,\hat{\xi},Z):= \sum_{j=1}^2C_{ij3l}(x,Z)\hat{\xi}_j~~~\text{   and   }~~~
   Q_{il}(x,\hat{\xi},Z):= \sum_{j,k=1}^2C_{ijkl}(x,Z)\hat{\xi}_j\hat{\xi}_k .
   \end{split}
\end{equation}
We follow \cite{tanuma} to define so-called limiting velocities. Let
$m$ and $n$ be orthogonal unit vectors in $\mathbb{R}^3$ which are
obtained by rotating the orthogonal unit vectors $\hat{\xi}$ and $e_3$
around their vector product $\hat{\xi} \times e_3$ by an angle
$\varrho\ (-\pi \leq \varrho \leq \pi)$:
\[m=m(\varrho)=\hat{\xi}\cos\varrho+e_3\sin\varrho,~~n=n(\varrho)=-\hat{\xi}\sin\varrho+e_3\cos\varrho.\]
For any $v>0$, we write
\[
   C_{ijkl}^v = C_{ijkl} - v^2 \hat{\xi}_j \hat{\xi}_k \delta_{il} .
\]
We let
\[
   T(\varrho,v,x,\hat{\xi},Z) = \sum_{j,k}C_{ijkl}^v(x,Z)n_jn_k ,
\]
\[
   R(\varrho;v,x,\hat{\xi},Z) = \sum_{j,k}C_{ijkl}^v(x,Z)m_j n_k\
~\text{and}~\
   Q(\varrho;v,x,\hat{\xi},Z) = \sum_{j,k}C_{ijkl}^v(x,Z)m_jm_k
\]
and note that
\[
   Q(\varrho+\tfrac{\pi}{2};\cdot) = T(\varrho;\cdot) ,\quad
   R(\varrho+\tfrac{\pi}{2};\cdot) = -R(\varrho;\cdot)^T ,\quad
   T(\varrho+\tfrac{\pi}{2};\cdot) = Q(\varrho;\cdot) .
\]

\begin{definition}\label{def:limvel}
The limiting velocity $v_L=v_L(x,\hat{\xi},Z)$ is the lowest velocity
for which the matrices $Q(\varrho;v,x,\hat{\xi},Z)$ and
$T(\varrho,v,x,\hat{\xi},Z)$ become singular for some angle $\varrho$:
\[
\begin{split}
	v_L(x,\hat{\xi},Z)=&\inf\{v>0\ :\ \exists \varrho :\ \det Q(\varrho;v,x,\hat{\xi},Z)=0\}\\
	=&\inf\{v>0\ :\ \exists \varrho :\ \det T(\varrho;v,x,\hat{\xi},Z)=0\}.
\end{split}
\]
\end{definition}

\begin{remark}
In the isotropic case, that is,
\[c_{ijkl}=\lambda \delta_{ij} \delta_{kl}
       + \mu (\delta_{ik} \delta_{jl} + \delta_{il} \delta_{jk}) ,
\]
where $\lambda, \mu$ are the two Lam\'{e} moduli, we have
\[
   C_{ijkl} = \hat{\lambda} \delta_{ij} \delta_{kl}
       + \hat{\mu} (\delta_{ik} \delta_{jl} + \delta_{il} \delta_{jk}) ,
\]
with $\hat{\mu}=\frac{\mu}{\rho}$ and $\hat{\lambda}=\frac{\lambda}{\rho}$.
The limiting velocity is given by
\[
   v_L(x,\hat{\xi},Z) = \sqrt{\hat{\mu}(x,Z)} .
\]
\end{remark}
\noindent
The spectrum of $H_0(x,\xi)$ is composed of a discrete spectrum
contained in $(0,v^2_L(x,\hat{\xi},Z_I)|\xi|^2)$ and an essential
spectrum, $[v^2_L(x,\hat{\xi},Z_I)|\xi|^2,\infty)$ (Figure
  \ref{fig:spectrumff}):

\begin{figure}[htp]
\centering
\includegraphics[width=3 in]{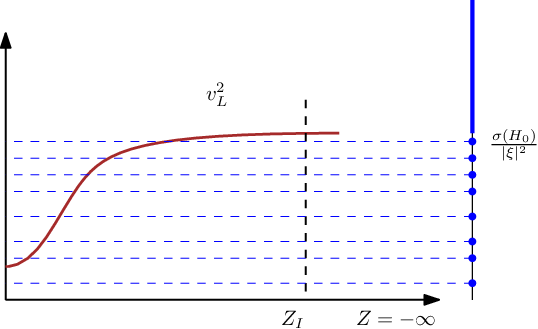}
\caption{Characterization of the spectrum}
\label{fig:spectrumff}
\end{figure}

\medskip\medskip

\begin{theorem}\label{existence}
Under Assumption $\ref{assu_tensor}$, the spectrum of $H_0(x,\xi)$ is
discrete below $v_L^2(x,\hat{\xi},Z_I)|\xi|^2$. The discrete spectrum
is nonempty if $|\xi|$ is sufficiently large.
\end{theorem}

\medskip\medskip

\begin{proof}
We need some notation and a few facts from operator theory for the
proof of this theorem. A brief review of those prerequisites are provided in
Appendix~B. For more details, we refer to \cite{RS}. We adopt the
classical Dirchlet-Neumann bracketing proof. Here we use $'$ to represent the partial derivative with respect to $Z$.  Consider
$h_0(x,\hat{\xi})$, and the corresponding quadratic form,
\begin{multline}
\mathfrak{h}_0(\varphi,\psi)=\int_{-\infty}^0 \sum_{i,l=1}^3\Big(h^2C_{i33l}\varphi'_i\overline{\psi'_l}-\ii h\sum_{j=1}^2C_{ij3l}\hat{\xi}_j\varphi'_i\overline{\psi_l}-\ii h\sum_{k=1}^2\hat{\xi}_k C_{i3kl}\varphi_i\overline{\psi'_l}\\
+\sum_{j,k=1}^2 C_{ijkl}\hat{\xi}_j\hat{\xi}_k\varphi_i\overline{\psi_l}\Big)\mathrm{d}Z
\end{multline}
on $H^1(\mathbb{R}^-)$, which is positive. We decompose $[Z_I,0]$ into $m$ intervals
$[Z_I,0]=\cup_{p=1}^m I_q$, the interiors, $(I_q)^\circ$, of which are
disjoint. We let $C_{ijkl}^{q,+}$ be a minimal element in the set,
\[\Pi^{q,+}=\{C_{ijkl} : C_{ijkl}\succeq C_{ijkl}(x,Z)\text{ for all }Z\in I_q\}.\] 
 (The order is defined as $C^1_{ijkl}\succeq C^2_{ijkl}$ if $\sum_{i,j,k,l=1}^3(C^1_{ijkl}-C^2_{ijkl})\epsilon_{ij}\epsilon_{kl}\geq 0$ for any symmetric matrix $\epsilon$. This is then a partially ordered set, and a minimal element exists by Zorn's lemma.)
Similarly, we define $C_{ijkl}^{q,-}$ to be a maximal element
of
\[\Pi^{q,-}=\{C_{ijkl}:C_{ijkl}\preceq C_{ijkl}(x,Z)\text{ for all }Z\in I_q\}.\]Furthermore, we introduce $C_{ijkl}^{m+1,\pm}=C_{ijkl}(x,Z_I)$ on $I_{m+1}=(-\infty, Z_I]$.
  
We consider
\begin{multline*}
\mathfrak{h}^{-,q}_0(\varphi,\psi)=\int_{I_q} \sum_{i,l=1}^3\Big(h^2C_{i33l}^{-,q}\varphi'_i\overline{\psi_l}'-\ii h\sum_{j=1}^2C_{ij3l}^{-,q}\hat{\xi}_j\varphi'_i\overline{\psi_l}-\ii h\sum_{k=1}^2\hat{\xi}_k C_{i3kl}^{-,q}\varphi_i\overline{\psi'_l}\\
+\sum_{j,k=1}^2 C_{ijkl}^{-,q}\hat{\xi}_j\hat{\xi}_k\varphi_i\overline{\psi_l}\Big)\mathrm{d}Z,
\end{multline*}
for $q=1,2,\cdots,m+1$, and let $h^{-,q}_0(x,\hat{\xi})$ be the unique
self-adjoint operator on $L^2(I_q)$ associated to the quadratic form
$\mathfrak{h}^{-,q}_0(\cdot,\cdot)$ (cf. Appendix B). This is equivalent to defining
\begin{multline}\label{operator}
(h^{-,q}_0(x,\hat{\xi}))_{il}=-h^2\frac{\partial}{\partial Z}C_{i33l}^{-,q}
\frac{\partial}{\partial Z}
- \ii h\sum_{j=1}^{2}
C_{ij3l}^{-,q}\hat{\xi}_j
\frac{\partial}{\partial Z}
- \ii h\sum_{k=1}^{2}
C_{i3kl}^{-,q} \frac{\partial }{\partial Z}
\xi_k\\
- \ii h\sum_{k=1}^{2}
\left( \frac{\partial}{\partial Z} C_{i3kl}^{-,q}\right)
\hat{\xi}_k 
+ \sum_{j,k=1}^{2}
C_{ijkl}^{-,q} \hat{\xi}_j \hat{\xi}_k ,
\end{multline}
on $\{u\in H^2(I_q):u~\text{with Neumann boundary condition}\}$.
Then
\[
   h_0(x,\hat{\xi}) \geq \bigoplus_{q=1}^{m+1} h^{-,q}_0(x,\hat{\xi}) .
\]
It follows that for each $q=1,\cdots, m$, $h^{-,q}_0(x,\hat{\xi})$ has
only a discrete spectrum, since it is a self-adjoint second-order
elliptic differential operator on a bounded interval.

The spectrum of $h^{-,m+1}_0(x,\hat{\xi})$ is also discrete below
$v^2_L(x,\hat{\xi},Z_I)$ \cite[Theorems 3.12, 3.13]{tanuma}.  We note
here that this property is related to the existence of surface waves
for a homogenous half space.  Thus the spectrum of $h_0(x,\hat{\xi})$
below $v^2_L(x,\hat{\xi},Z_I)$ must be discrete (see the Lemma on
pp. 270, Vol. IV of \cite{RS}.)

We define $h^{+,q}_0(x,\hat{\xi})$ as the unique self-adjoint operator
on $L^2(I_q)$, the quadratic form of which is the closure of the form
\begin{multline*}
\mathfrak{h}^{+,q}_0(\varphi,\psi)=\int_{I_q} \sum_{i,l=1}^3\Big(h^2C_{i33l}^{+,q}\varphi'_i\overline{\psi_l}'-\ii h\sum_{j=1}^2C_{ij3l}^{+,q}\hat{\xi}_j\varphi'_i\overline{\psi_l}-\ii h\sum_{k=1}^2\hat{\xi}_k C_{i3kl}^{+,q}\varphi_i\overline{\psi'_l}\\
+\sum_{j,k=1}^2 C_{ijkl}^{+,q}\hat{\xi}_j\hat{\xi}_k\varphi_i\overline{\psi_l}\Big)\mathrm{d}Z,
\end{multline*}
with domain $C_0^\infty(I_q)$. This is equivalent to defining
\begin{multline}\label{operator1}
(h_0^{+,q}(x,\hat{\xi}))_{il}=-h^2\frac{\partial}{\partial Z}C_{i33l}^{+,q}
\frac{\partial}{\partial Z}
- \ii h\sum_{j=1}^{2}
C_{ij3l}^{+,q}\hat{\xi}_j
\frac{\partial}{\partial Z}
- \ii h\sum_{k=1}^{2}
C_{i3kl}^{+,q} \frac{\partial }{\partial Z}
\hat{\xi}_k\\
- \ii h\sum_{k=1}^{2}
\left( \frac{\partial}{\partial Z} C_{i3kl}^{+,q}\right)
\hat{\xi}_k 
+ \sum_{j,k=1}^{2}
C_{ijkl}^{+,q} \hat{\xi}_j \hat{\xi}_k ,
\end{multline}
on $H^1_0(I_q)\cap H^2(I_q)$. Then
\[h_0(x,\hat{\xi})\leq \bigoplus_{q=1}^{m+1} h^{+,q}_0(x,\hat{\xi}).\]
$h^{+,j}_0(x,\hat{\xi})$, for each $j=1,\cdots,m$, has a discrete and
real spectrum. The operator $h^{+,m+1}_0(x,\hat{\xi})$ also has a
discrete spectrum below $v_L(x,\hat{\xi},Z_I)$, since
$h^{+,m+1}_0(x,\hat{\xi})\geq h^{-,m+1}_0(x,\hat{\xi})$.

We suppose now that the decomposition is fine enough, and that
$C_{ijkl}^{+,q}$ for some $q$ has a limiting velocity
$v_L^{+,q}(x,\hat{\xi}) < v_L(x,\hat{\xi},Z_I)$. (There exists $Z^*
\in (Z_I,0]$ such that
  $v_L(x,\hat{\xi},Z^*)<v_L(x,\hat{\xi},Z_I)$. Suppose $Z^* \in I_q$
  for some $q$. For any $\varepsilon>0$, if $|I_q|$ is small enough,
  $C_{ijkl}(x,Z^*)+\varepsilon
  (\delta_{ik}\delta_{jl}+\delta_{ij}\delta_{kl}+\delta_{il}\delta_{kj})
  \succeq C_{ijkl}(x,Z),\forall Z\in I_q$. This show nonemptyness of
  $\Pi^{q,+}$. We can take $\varepsilon$ small enough, so that
  $C_{ijkl}^{+,q}=C_{ijkl}(x,Z^*)+\varepsilon
  (\delta_{ik}\delta_{jl}+\delta_{ij}\delta_{kl}+\delta_{il}\delta_{kj})$
  has limiting velocity
  $v_L^{+,q}(x,\hat{\xi})<v_L(x,\hat{\xi},Z_I)$.)
We let $\lambda_{q,1}$ be the first eigenvalue of $h^{+,q}_0$, and
claim that
$\lambda_{q,1}<(v^{+,q}_L(x,\hat{\xi}))^2+\mathcal{O}(h)$. We prove
this claim next.

We let $I_q = [a,b]$. We note that
\[
   h^{+,q}_0 u = ((v_L^{+,q}(x,\hat{\xi}))^2+h) u
\]
has a bounded solution, $u = \mathbf{a} e^{\frac{\ii p Z}{h}}$ with $\mathbf{a}\in\mathbb{C}^3$, $p\in\mathbb{R}$. To
construct a solution of this form, we carry out a substitution and
obtain
\begin{equation}\label{det1a}
   [T^{+,q} p^2 + (R^{+,q} + (R^{+,q})^T) p
            + Q^{+,q} - v_L^{+,q}(x,\hat{\xi})^2-h] \mathbf{a} = 0 .
\end{equation}
Here, $T^{+,q},R^{+,q},Q^{+,q}$ are defined as in (\ref{TRQ}) for
$C_{ijkl}^{q,+}$. The above system has a non-trivial solution if and
only if
\begin{equation}
\label{deta}
   \det[T^{+,q} p^2 + (R^{+,q} + (R^{+,q})^T) p + Q^{+,q}
              - v_L^{+,q}(x,\hat{\xi})^2-h] = 0 .
\end{equation}
Since $v_L^{+,q}(x,\hat{\xi}))^2+h > v_L^{+,q}(x,\hat{\xi}))^2$, there
exists a real-valued solution $p$ to (\ref{deta}) (\cite{tanuma}, Lemma 3.2).
.

We let $v=u$ on $[a+h,b-h]$ and $v \in C^\infty_0(I_q)$;
$v$ can be constructed such that $\|v\|_{L^2([a,a+h]\cup[b-h,b])}$ $\leq
Ch^{1/2}$, $\|v'\|_{L^2([a,a+h]\cup[b-h,b])}\leq Ch^{-1/2}$,
$\|v''\|_{L^2([a,a+h]\cup[b-h,b])}\leq Ch^{-3/2}$. Then
$C_1<\|v\|<C_2$ for some $C_1,C_2>0$ independent of $h$. After
renormalization, we can assume that $\|v\|_2 = 1$. Then
\begin{multline}
   ( h^{+,q}_0 v, v)
   \le \int_{a+h}^{a-h} h_0^{+,q} v \cdot \overline{v}\mathrm{d}Z
       + \left( \int_{a}^{a+h}+\int_{b-h}^{b} \right) h_0^{+,q}
                v \cdot \overline{v}\mathrm{d}Z
\\
   \le (v_L^2(x,\hat{\xi},Z^*) + h) \|v\|^2_{L^2([a+h,b-h])}
\\[0.25cm]
       + C h^2 \|v\|_{L^2([a,a+h]\cup[b-h,b])}
                  \|v''\|_{L^2([a,a+h]\cup[b-h,b])}
\\[0.25cm]
   + C h \|v\|_{L^2([a,a+h]\cup[b-h,b])}
        \|v'\|_{L^2([a,a+h]\cup[b-h,b])}
   + C \|v\|_{L^2([a,a+h]\cup[b-h,b])}
                  \|v\|_{L^2([a,a+h]\cup[b-h,b])}
\\[0.25cm]
   \le (v_L^2(x,\hat{\xi},Z^*) + h) (1 - C h) + C h .
\end{multline}
When $h$ is small, that is, $|\xi|$ is large, $h^{+,q}_0$ has
eigenvalues below $v^2_L(x,\hat{\xi},Z_I)$. Then the discrete spectrum
of $h_0(x,\hat{\xi})$ is nonempty.
\end{proof}\\

\begin{theorem}\label{ess}
The essential spectrum of $H_0(x,\xi)$ is given by
$[v_L^2(x,\hat{\xi},Z_I) |\xi|^2,+\infty)$.
\end{theorem}\\

\begin{proof}
We let $H_L(x,\xi)$ be the operator given by
\begin{multline}\label{operatorp}
(H_L(x,\xi))_{il}=-h^2\frac{\partial}{\partial Z}C_{i33l}(x,Z_I)
\frac{\partial}{\partial Z}
- \ii h\sum_{j=1}^{2}
C_{ij3l}(x,Z_I)\hat{\xi}_j
\frac{\partial}{\partial Z}
- \ii h\sum_{k=1}^{2}
C_{i3kl}(x,Z_I) \frac{\partial }{\partial Z}
\xi_k\\
- \ii h\sum_{k=1}^{2}
\left( \frac{\partial}{\partial Z} C_{i3kl}(x,Z_I)\right)
\hat{\xi}_k 
+ \sum_{j,k=1}^{2}
C_{ijkl}(x,Z_I) \hat{\xi}_j \hat{\xi}_k ,
\end{multline}
with (locally) constant coefficients $C_{ijkl}(x,Z_I)$ defined on $L^2(\mathbb{R}^-)$ supplemented with the Neumann boundary
condition. This operator is self adjoint in $L^2(\mathbb{R}^-)$. We
note that $H_0(x,\xi) - H_L(x,\xi)$ is symmetric and is a relatively
compact perturbation of $(H_L(x,\xi))^2$. This follows from the
observation that
\[
   \left(H_0(x,\xi) - H_L(x,\xi)\right) (H_L(x,\xi)^2 + \ii)^{-1}
\]
maps $L^2(\mathbb{R^-})$ to the set
\[
   \mathcal{L} = 
     \{f\ :\ f \in H^2(\mathbb{R}^-),
                         f\text{ is supported in }[Z_I,0]\} ,
\]
while any bounded subset of $\mathcal{L}$ is compact in
$L^2(\mathbb{R}^-)$. Thus $H_0(x,\xi)$ and $H_L(x,\xi)$ share the same
essential spectrum \cite[Corollary 3 in XIII.4]{RS}.

We claim that $H_L(x,\xi)$ has $[v^2_L(x,\hat{\xi},Z_I)
  |\xi|^2,\infty)$ in its essential spectrum. To show this, we
  construct a singular sequence $\{v_k\}_{k=1}^\infty$ satisfying Theorem \ref{eigcri} for any
  $\Lambda > v^2_L(x,\hat{\xi},Z_I) |\xi|^2$. In this proof, we use
  the simplified notation, $T = T(x,Z_I)$, $R = R(x,\hat{\xi},Z_I)$ and
  $Q = Q(x,\hat{\xi},Z_I)$. Also, we use $\Lambda = v^2 |\xi|^2$. We
  first seek solutions, ignoring boundary conditions, of the form
\[
   u = \mathbf{a} \, e^{\frac{\ii p Z}{h}} .
\]

Substituting this form into the equation, we obtain
\begin{equation}\label{det1}
   [T p^2 + (R + R^T) p + Q - v^2] \mathbf{a} = 0 .
\end{equation}
The above system has a non-trivial solution if and only if 
\begin{equation}
\label{det}
   \det[T p^2 + (R + R^T) p + Q -v^2] = 0.
\end{equation}
When $v^2>v^2_L(x,\hat{\xi},Z_I)$, there exists a real-valued solution
$p$ to (\ref{det}) (\cite{tanuma}, Lemma 3.2).

We take $u_k = \mathbf{a} e^{\frac{\ii p Z}{h}} \phi_k(Z)$, with
$\phi_k \in C_0^\infty(\mathbb{R}^-)$, $\supp\phi_k(Z)\subset [1,k]$,
$\phi_k(Z) = 1$ on $[2,k-1]$ and $|\partial^\alpha\phi_k| \leq
C,|\alpha| \leq 2$. Then $\|u_k\|_2 \geq C k$ for some constant $C >
0$. We note that $u_k \in \mathcal{D}$. With $v_k = \frac{u_k}{\|u_k\|}$ we find that
\[
   \| H_L(x,\xi) v_k - v^2 |\xi|^2 v_k\|
                    \leq \frac{C}{k} \rightarrow 0 .
\]
Thus $v^2|\xi|^2$ is in the spectrum of $H_0(x,\xi)$.

This shows that any $\Lambda > v_L^2(x,\hat{\xi},Z_I)|\xi|^2$ is in
the essential spectrum of $H_L(x,\xi)$. Then, by Theorem \ref{existence}, $[v_L^2(x,\hat{\xi},Z_I) |\xi|^2,\infty)$ is the essential spectrum
  of $H_0(x,\xi)$.
\end{proof}

\medskip\medskip

\noindent
We emphasize that, here, the essential spectrum is not necessarily a
continuous spectrum. It may contain eigenvalues.

\section{Surface-wave modes associated with exponentially
         decaying eigenfunctions}
\label{sec:4}
\begin{figure}[h]
\centering
\includegraphics[width=5 in]{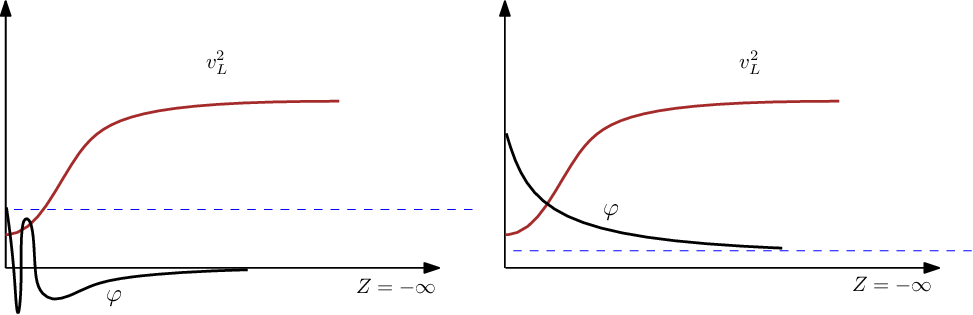}
\caption{Behaviors of eigenfunctions. The dashed line indicates the
  eigenvalue }
\label{eigenf}
\end{figure}

~\\
When the limiting velocity is minimal at the surface, that is,
\begin{equation}\label{BarL}
   v_L(x,\hat{\xi},0) = \inf_{Z\leq 0} v_L(x,\hat{\xi},Z) ,
\end{equation}
there may still exist an eigenvalue of $H_0(x,\xi)$ below
$v^2_L(x,\hat{\xi},0) |\xi|^2$. The corresponding eigenfunction
exponentially decays ``immediately'' from the surface. This is
somewhat different from the behavior of the eigenfunctions with
eigenvalues above $v^2_L(x,\hat{\xi},0)|\xi|^2$ (See Figure
\ref{eigenf}). Here, we discuss the existence of such an eigenvalue.

For the principal symbol (cf.~(\ref{principle})) of
$h_0(x,\hat{\xi})$, for $0<v^2<v_L^2(x,\hat{\xi},0)$, we have
\[
   \det[T \zeta^2 + (R+R^T) \zeta + Q - v^2] = 0 ,
\]
and viewed as a polynomial in $\zeta$ has $6$ non-real roots that
appear in conjugate pairs. Suppose that $\zeta_1,\zeta_2,\zeta_3$ are
three roots with negative imaginary parts, and
$\gamma\subset\mathbb{C}_-$ is a continuous curve enclosing
$\zeta_1,\zeta_2,\zeta_3$. We let
\begin{equation}\label{D1}
   S_1(x,Z,v,\hat{\xi}))
          = T^{-1/2}\tilde{S}_1(x,Z,v,\hat{\xi})T^{1/2} ,
\end{equation}
where
\begin{equation}\label{D2}
   \tilde{S}_1(x,Z,v,\hat{\xi}) :=
      \left(\oint_\gamma\zeta \tilde{M}(x,Z,v,\hat{\xi},\zeta)^{-1}
            \mathrm{d}\zeta\right)
      \left(\oint_\gamma \tilde{M}(x,Z,v,\hat{\xi},\zeta)^{-1}
            \mathrm{d}\zeta\right)^{-1} ,
\end{equation}
and
\[
   \tilde{M}(x,Z,v,\hat{\xi},\zeta)
      = T^{-1/2}[T\zeta^2+(R+R^T)\zeta+Q-v^2I]T^{-1/2} .
\]
Then we have \cite{Hansen,tanuma}
\[
   T \zeta^2 + (R+R^T) \zeta + Q - v^2 I
      = (\zeta-S_1^*(x,Z,v,\hat{\xi})) T
                           (\zeta-S_1(x,Z,v,\hat{\xi})) ,
\]
with $\mathrm{Spec}(S_1) \subset \mathbb{C}_-$, for
$0<v^2<v_L^2(x,\hat{\xi},0)$.

We define an operator,
\begin{equation}\label{unpert}
   l_0(x,\hat{\xi},v)(Z,hD_Z) = (hD_Z-S_1^*)T(hD_Z-S_1) .
\end{equation}
A solution to $l_0(x,\hat{\xi},v)(Z,hD_Z) \varphi = 0$ in
$H^2(\mathbb{R}^-)$ needs to satisfy
\begin{equation}\label{reduced}
   (hD_Z-S_1) \varphi = 0 .
\end{equation}
This follows from letting $f = T (hD_Z-S_1) \varphi$, and noting that
$(hD_Z-S_1^*) f = 0$. If $(hD_Z-S_1) \varphi = T^{-1}f(Z_0) \neq 0$
for some $Z_0 \in\mathbb{R}^-$, then
\[
   f(Z) = f(Z_0) e^{\ii\frac{ S_1^*(Z-Z_0)}{h}}
\]
is exponentially growing as $Z \rightarrow -\infty$, since
$\mathrm{Spec}(S_1^*)\subset\mathbb{C}_+$. This contradicts the fact
that $\varphi \in H^2(\mathbb{R}^-)$.

We now introduce the surface impedance tensor,
\begin{equation}\label{ZZ}
   \mathcal{Z}(x,v,\hat{\xi})
      = -\ii \, (T(x,0,\hat{\xi})S_1(x,0,v,\hat{\xi})
              + R^T(x,0,\hat{\xi})) .
\end{equation}
The basic properties of this tensor are summarized in

\medskip\medskip

\begin{proposition}[\cite{Hansen}]\label{Hansen}
For $0 \leq v < v_L(x,\hat{\xi},0))$ the following holds true
\begin{enumerate}[label=(\roman*)]
\item $\mathcal{Z}(x,v,\hat{\xi})$ is Hermitian.
\item $\mathcal{Z}(x,0,\hat{\xi})$ is positive definite.
\item The real part of $\mathcal{Z}(x,v,\hat{\xi})$ is positive
  definite.
\item At most one eigenvalue of $\mathcal{Z}(x,v,\hat{\xi})$ is
  non-positive.
\item The derivative,
  $\frac{\mathrm{d}\mathcal{Z}(x,v,\hat{\xi})}{\mathrm{d}v}$, is
  negative definite.
\item The limit $\mathcal{Z}(v_L(x,\hat{\xi},0))) = \lim_{v\uparrow
  v_L(x,\hat{\xi},0))}\mathcal{Z}(x,v,\hat{\xi})$ exists.
\end{enumerate}
\end{proposition}

\medskip

To make this paper more self-contained, we include a proof of the
above proposition. The properties of $\mathcal{Z}$ can be studied also
by Stroh's formalism and are important for the existence of subsonic
Rayleigh wave in a homogeneous half-space; see \cite{tanuma} for a
comprehensive overview.

\begin{proof}
We use the Riccati equation
\begin{equation}\label{Ri}
(\mathcal{Z}(x,v,\hat{\xi})-\ii R(x,\hat{\xi},0))T(x,0,\hat{\xi})^{-1}(\mathcal{Z}(x,v,\hat{\xi})+\ii R(x,\hat{\xi},0)^T)=Q(x,0,\hat{\xi})-v^2I
\end{equation}
as in \cite{MF,Hansen}. We note that
\[S_1(x,0,v,\hat{\xi})=\ii T(x,0,\hat{\xi})^{-1}(\mathcal{Z}(x,v,\hat{\xi})+\ii R^T).\]
Taking the adjoint of $(\ref{Ri})$, we find that
\begin{equation}\label{Ri1}
(\mathcal{Z}(x,v,\hat{\xi})^*-\ii R(x,\hat{\xi},0))T(x,0,\hat{\xi})^{-1}(\mathcal{Z}(x,v,\hat{\xi})^*+\ii R(x,\hat{\xi},0)^T)=Q(x,0,\hat{\xi})-v^2I.
\end{equation}
Subtracting $(\ref{Ri1})$ from $(\ref{Ri})$, we obtain
\[S_1(x,0,v,\hat{\xi})^*(\mathcal{Z}(x,v,\hat{\xi})-\mathcal{Z}(x,v,\hat{\xi})^*)-(\mathcal{Z}(x,v,\hat{\xi})-\mathcal{Z}(x,v,\hat{\xi})^*)S_1(x,0,v,\hat{\xi})=0.\]
The above Sylvester equation is nonsingular, because $S_1(x,0,v,\hat{\xi})$ and $S_1(x,0,v,\hat{\xi})^*$ have disjoint spectra. Therefore
\[\mathcal{Z}(x,v,\hat{\xi})=\mathcal{Z}(x,v,\hat{\xi})^*.\]
This proves (i).

Differentiating  $(\ref{Ri})$ in $v$, we get
\[\ii S_1(x,0,v,\hat{\xi})^*\frac{\mathrm{d}}{\mathrm{d}v}\mathcal{Z}(x,v,\hat{\xi})-\ii \frac{\mathrm{d}}{\mathrm{d}v}\mathcal{Z}(x,v,\hat{\xi})S_1(x,0,v,\hat{\xi})=-2v I.\]
The above Lyapunov-Sylvester equation has the unique solution
\[\frac{\mathrm{d}}{\mathrm{d}v}\mathcal{Z}(x,v,\hat{\xi})=-2v\int_0^\infty\exp{(-\ii t S_1(x,0,v,\hat{\xi})^*)}\exp{(\ii t S_1(x,0,v,\hat{\xi}))}\mathrm{d}t.\]
It is clear that
$\frac{\mathrm{d}}{\mathrm{d}v}\mathcal{Z}(x,v,\hat{\xi})$ is negative definite since $S_1(x,0,v,\hat{\xi})\subset \mathbb{C}_-$. This proves (v).

In order to prove (ii), we assume that $w=\left(\begin{array}{c}w_1\\w_2\\w_3\end{array}\right)$ is a complex vector such that
\[w^T\mathcal{Z}(x,0,\hat{\xi})\overline{w}\leq 0.\]
Let
\[u(x,\hat{\xi},Z)=e^{\ii\frac{ZS_1(x,0,0,\hat{\xi})}{h}}w,~~\text{on}~Z\in(-\infty,0].\]
Then
\[-\left(hT(x,0,\hat{\xi})\frac{\partial}{\partial Z}+\ii R(x,0,\hat{\xi})^T\right)u\vert_{Z=0}=\mathcal{Z}(x,0,\hat{\xi})w.\]
We have the energy identity
\begin{multline*}
\int_{-\infty}^0 \sum_{i,l=1}^3\Big(h^2C_{i33l}(x,0)u'_i\overline{u'_l}-\ii h\sum_{j=1}^2C_{ij3l}(x,0)\hat{\xi}_ju'_i\overline{u_l}-\ii h\sum_{k=1}^2\hat{\xi}_k C_{i3kl}(x,0)u_i\overline{u'_l}\\
+\sum_{j,k=1}^2 C_{ijkl}(x,0)\hat{\xi}_j\hat{\xi}_ku_i\overline{u_l}\Big)\mathrm{d}Z=w^T\mathcal{Z}(x,0,\hat{\xi})\overline{w}\leq 0.
\end{multline*}
By positivity of the differential operator, $u=0$. Thus $w=0$, and we have proved (ii).

With (\ref{D1}), (\ref{D2}) and (\ref{ZZ}), we obtain
\[
\ii\mathcal{Z}(x,v,\hat{\xi})\left(\oint_{\gamma_R} M(x,0,v,\hat{\xi},\zeta)^{-1}
            \mathrm{d}\zeta\right)=\oint_{\gamma_R}\left(\zeta T(x,0,\hat{\xi})+R^T(x,0,\hat{\xi})\right) M(x,0,v,\hat{\xi},\zeta)^{-1}
            \mathrm{d}\zeta,
\]
where $\gamma_R$ is the closed contour consisting of $[-R,R]$ on the real line and the arc $Re^{\ii\theta}$, $\theta\in[0,\pi]$, for $R$ large enough. We note that 
\[T(x,0,\hat{\xi})M(x,0,v,\hat{\xi},\zeta)^{-1}=\zeta^{-2}I+\mathcal{O}(|\zeta|^{-3}).\]
We then have, by letting $R\rightarrow\infty$,
\[\ii
\mathcal{Z}(x,v,\hat{\xi})\left(\int_{-\infty}^{+\infty} M(x,0,v,\hat{\xi},\zeta)^{-1}
            \mathrm{d}\zeta\right)=\ii\pi I+\int_{-\infty}^{+\infty}\left(\zeta T(x,0,\hat{\xi})+R^T(x,0,\hat{\xi})\right) M(x,0,v,\hat{\xi},\zeta)^{-1}
            \mathrm{d}\zeta.
\]
The integral on the right-hand side is a real-valued matrix, and
\[
\int_{-\infty}^{+\infty} M(x,0,v,\hat{\xi},\zeta)^{-1}
            \mathrm{d}\zeta
\]
is a positive definite matrix. Thus $\Re Z(x,v,\hat{\xi})$ is positive
definite. This proves (iii).

We prove (iv) by contradiction: Assume that
$\mathcal{Z}(x,v,\hat{\xi})$ has only one positive eigenvalue for some
$v$. The corresponding eigenspace would be one-dimensional. We can
then choose a real-valued vector $w\neq 0$ orthogonal to the
eigenspace such that $w^T\mathcal{Z}(x,v,\hat{\xi}) w\leq 0$. This
contradicts the fact that $\Re \mathcal{Z}(x,v,\hat{\xi})$ is positive
definite, and proves (iv).

Finally, we prove (vi). Let $\|\mathcal{Z}(x,v,\hat{\xi})\|$ be the operator norm of $\mathcal{Z}(x,v,\hat{\xi})$. Then $\|\mathcal{Z}(x,v,\hat{\xi})\|$ is equal to the maximum of the absolute values of eigenvalues of $\mathcal{Z}(x,v,\hat{\xi})$. By (v), we have
\[\mathcal{Z}(x,v,\hat{\xi})\preceq  \mathcal{Z}(x,0,\hat{\xi}).\]
Here $A \preceq B$, with Hermitian matrices $A$ and $B$, means
\[w^T(A-B)\overline{w}\leq 0\]
for any complex vector $w$. Thus the eigenvalues of $\mathcal{Z}(x,v,\hat{\xi})$ are not greater than $\|\mathcal{Z}(x,0,\hat{\xi})\|$ for $0\leq v\leq v_L(x,\hat{\xi},0))$. The sum of eigenvalues of $\mathcal{Z}(x,v,\hat{\xi})$ is positive, since it is equal to $\text{trace}(\mathcal{Z}(x,v,\hat{\xi}))$, which is in turn equal to  $\text{trace}(\Re\mathcal{Z}(x,v,\hat{\xi}))$. Therefore, the possibly negative eigenvalue of $\mathcal{Z}(x,v,\hat{\xi})$ is less than the sum of the two positive eigenvalues. Then
\[\|\mathcal{Z}(x,v,\hat{\xi})\|\leq 2\|\mathcal{Z}(x,0,\hat{\xi})\|
\]
and we can take the limit in (vi).
\end{proof}

We now introduce the generalized Barnett-Lothe condition
\cite{Nakamura}: $(x,\xi)$ satisfies the generalized Barnett-Lothe
condtion if
\[
   \lim_{v\uparrow v_L(x,\hat{\xi},0)}
            \det \mathcal{Z}(x,v,\hat{\xi}) < 0 ,
\]
or
\[
   \lim_{v\uparrow v_L(x,\hat{\xi},0)}
     [(\mathrm{trace}\mathcal{Z}(x,v,\hat{\xi}))^2
           -\mathrm{trace}\mathcal{Z}^2(x,v,\hat{\xi})] < 0 .
\]
We note that in the isotropic case, this condition is satisfied for
all $(x,\hat{\xi})$.

We let $\varphi$ be a solution of (\ref{reduced}). Then the Neumann
boundary value of $\varphi$ is given by
\[
\begin{split}
   \ii \, (-hT(x,0)D_Z - R(x,0,\hat{\xi})^T)\varphi_0(0)
   &=-\ii\, (T(x,0)S_1(x,0,v,\hat{\xi})
                 + R(x,0,\hat{\xi})^T) \varphi_0(0)
\\
   &=\mathcal{Z}(x,v,\hat{\xi})\varphi_0(0) .
\end{split}
\]
The generalized Barnett-Lothe condition guarantees that
$\mathcal{Z}(x,v_L(x,\hat{\xi},0),\hat{\xi})$ has one negative
eigenvalue. Therefore, with properties of $\mathcal{Z}(x,v,\hat{\xi})$
established in Proposition \ref{Hansen}, we have

\medskip\medskip

\begin{proposition}\label{R_speed}
There is a unique $v_0$ with $0 \leq v_0^2<v_L^2(x,\hat{\xi},0)$ such
that $\det \mathcal{Z}(x,v_0,\hat{\xi}) = 0$, if the generalized
Barnett-Lothe condition holds.
\end{proposition}

\medskip\medskip

At the $v_0$ given in the proposition above, (\ref{reduced}) has a
solution $\varphi_0$ that satisfies Neumann boundary condtion, that
is, $\varphi_0 \in \mathcal{D}$. We further assume
$\|\varphi_0\|_2=1$. We emphasize, here, that the existence of such a
$v_0$ depends on the stiffness tensor at boundary only. This locality
behavior was first observed by Petrowsky \cite{Petrowsky}.

We find that
\[
   (h_0(x,\hat{\xi})(Z,hD_Z) - v^2_0 - l_0(x,\hat{\xi},v_0)(Z,hD_Z))
           \varphi_0 = h (D_Z (TS_1+R^T)) \varphi_0
\]
by straightforward calculations, noting that
\[
   h_0(x,\hat{\xi})(Z,hD_Z)\varphi_0
      =h^2 D_Z(TD_Z\varphi_0)+h(R(D_Z\varphi_0
                      + D_Z (R^T \varphi_0)) + Q \varphi_0 ,
\]
\[
   S_1^*TS_1=Q-v_0^2
\]
and
\[
   -S_1^*T-TS_1=R+R^T .
\]
Then
\[
   \mathfrak{h}_0(x,\hat{\xi})(\varphi_0,\varphi_0)
      = (h_0(x,\hat{\xi})(Z,hD_Z)\varphi_0,\varphi_0) \leq v^2_0 + C h ,
\]
since $\varphi_0$ satisfies the Neumann boundary condition.  Thus the
first eigenvalue $v^2_1(x,\xi)$ of \\ $h_0(x,\hat{\xi})(Z,hD_Z)$
satisfies
\[
   v^2_1(x,\xi) = \min_{\|\varphi\|=1}
       \mathfrak{h}_0(x,\hat{\xi})(\varphi,\varphi) \leq v^2_0 + C h .
\]
Therefore,

\medskip\medskip

\begin{theorem}
Assume that $(\ref{BarL})$ holds. If the generalized Barnett-Lothe
condtion holds, and $|\xi|$ is sufficiently large, then there exists
an eigenvalue $\Lambda_1 < v_L^2(x,\hat{\xi},0)|\xi|^2$ of
$H_0(x,\xi)$.
\end{theorem}

\medskip\medskip

This particular eigenvalue may exist even when $C_{ijkl}(x,\cdot)$ is
constant (independent of $x$). It was first found by Rayleigh himself
for the isotropic case.  The uniqueness of such subsonic eigenvalue
can be guaranteed if we further assume that $C_{ijkl}(x,Z)$ is
non-increasing as $Z$ with respect to the order defined in the proof
of Theorem \ref{existence}.

\section{Decoupling into Rayleigh and Love waves}
\label{sec:5}
\subsection{Isotropic case}

For the case of isotropy,
\begin{multline}
H_{0} u =
\left[\left(\begin{array}{ccc}
-\pdpd{}{Z} (\hat{\mu}(x,Z)\pdpd{}{Z} )&
0 &
0 \\
0 &
-\pdpd{}{Z} (\hat{\mu}(x,Z) \pdpd{}{Z}) &
0 \\
0 &
0 &
-\pdpd{}{Z}
((\hat{\lambda} + 2\hat{\mu})(x,Z) \pdpd{}{Z})
\end{array}\right)\right.
\\
- \ii \left(\begin{array}{ccc}
0 &
0 &
\hat{\lambda}(x,Z) \xi_1
\pdpd{}{Z} \\
0 &
0 &
\hat{\lambda}(x,Z)\xi_2
\pdpd{}{Z} \\
\hat{\mu}(x,Z)\xi_1\pdpd{}{Z} &
\hat{\mu}(x,Z)\xi_2\pdpd{}{Z} &
0
\end{array}\right)
\\
- \ii  \left(\begin{array}{ccc}
0 &
0 &
\hat{\mu}(x,Z) \xi_1
\pdpd{}{Z} \\
0 &
0 &
\hat{\mu}(x,Z)\xi_2
\pdpd{}{Z} \\
\hat{\lambda}(x,Z)\xi_1\pdpd{}{Z} &
\hat{\lambda}(x,Z)\xi_2\pdpd{}{Z} &
0
\end{array}\right)
\\
- \ii \left(\begin{array}{ccc}
0 &
0 &
\left(\pdpd{}{Z} \hat{\mu}(x,Z)\right)
\xi_1  \\
0 &
0 &
\left(\pdpd{}{Z} \hat{\mu}(x,Z)\right)
\xi_2  \\
\left(\pdpd{}{Z} \hat{\lambda}(x,Z)\right)
\xi_1 &
\left(\pdpd{}{Z} \hat{\lambda}(x,Z)\right)
\xi_2  &
0
\end{array}\right)
\\
\left.
+ \left(\begin{array}{ccc}
 (\hat{\lambda} + 2\hat{\mu}) \, \xi_1^2
+ \hat{\mu}\, \xi_2^2  \, &
(\hat{\lambda} + \hat{\mu}) \, \xi_1 \xi_2 \, &
0 \\
(\hat{\lambda} + \hat{\mu}) \, \xi_1 \xi_2 \,  &
 (\hat{\lambda} + 2\hat{\mu})\, \xi_2^2
+ \hat{\mu} \, \xi_1^2\, &
0 \\
0 &
0 &
\hat{\mu} \, (\xi_1^2 + \xi_2^2) \, 
\end{array}\right)\right]
\left(\begin{array}{c}u_1\\ u_2\\u_3\end{array}\right).
\end{multline}
We introduce an orthogonal matrix
\[P(\xi)=\left(\begin{array}{ccc}
|\xi|^{-1}\xi_1 & |\xi|^{-1}\xi_2 & 0\\
|\xi|^{-1}\xi_2 &-|\xi|^{-1}\xi_1 &0 \\
 0 & 0 & 1
\end{array}\right)\]
Using the substitution $u=P(\xi)\varphi$, with
$\varphi=(\varphi_1,\varphi_2,\varphi_3)$, we obtain the equations for
the eigenfunctions,
\begin{equation}\label{love}
-\pdpd{}{Z} \hat{\mu} \pdpd{\vphi_2}{Z}
+ \hat{\mu} \, |\xi|^2 \vphi_2
= \Lambda\vphi_2 ,
\end{equation}
\begin{equation}\label{love_b}
\pdpd{\vphi_2}{Z}(0) = 0 ,
\end{equation}
for Love waves, and
\begin{equation}\label{Rayleigh1}
-\frac{\partial}{\partial Z}\hat{\mu}\frac{\partial\varphi_1}{\partial Z}-\ii|\xi|\left(\frac{\partial}{\partial Z}(\hat{\mu}\varphi_3)+\hat{\lambda}\frac{\partial}{\partial Z}\varphi_3\right)+(\hat{\lambda}+2\hat{\mu})|\xi|^2\varphi_1=\Lambda\varphi_1 ,
\end{equation}
\begin{equation}\label{Rayleigh2}
-\frac{\partial}{\partial Z}(\hat{\lambda}+2\hat{\mu})\frac{\partial\varphi_3}{\partial Z}-\ii|\xi|\left(\frac{\partial}{\partial Z}(\hat{\lambda}\varphi_1)+\hat{\mu}\frac{\partial}{\partial Z}\varphi_1\right)+\hat{\mu}|\xi|^2\varphi_3=\Lambda\varphi_3 ,
\end{equation}

\begin{eqnarray}
\ii|\xi| \vphi_3(0) + \pdpd{\vphi_1}{Z}(0) &=& 0 ,\label{Rayleighboundary1}
\\
\ii\hat{\lambda} |\xi| \vphi_1(0) + (\hat{\lambda} + 2\hat{\mu}) \pdpd{\vphi_3}{Z}(0)\label{Rayleighboundary2}
&=& 0 ,
\end{eqnarray}
for Rayleigh waves. So we have two types of surface-wave modes, Love
and Rayleigh waves, decoupled up to principal parts in an isotropic
medium. The lower-order terms can be constructed following the proof
of Theorem \ref{conjugation_sys} leading to coupling.

\subsection{Transversely isotropic case}

We consider a transversely isotropic medium having the surface normal
direction as symmetry axis. Then the nonzero components of $C$ are
\[C_{1111},~C_{2222},~C_{3333},~C_{1122},~C_{1133},~C_{2233},~C_{2323},~C_{1313},~C_{1212}\]
and
\[C_{1111}=C_{2222}, ~~C_{1133}=C_{2233},\]
\[C_{2323}=C_{1313},~~C_{1212}=\frac{1}{2}(C_{1111}-C_{1122}).\]
Using the substitution $u=P(\xi)\varphi$, again,
we obtain the equations for the eigenfunctions from $P^{-1}(\xi)H_0(x,\xi)P(\xi)\varphi =\Lambda \varphi$:
\begin{equation}\label{love_ti}
-\pdpd{}{Z} C_{1313} \pdpd{\vphi_2}{Z}
+ C_{1212} \, |\xi|^2 \vphi_2
= \Lambda\vphi_2 ,
\end{equation}
\begin{equation}\label{love_b_ti}
\pdpd{\vphi_2}{Z}(0) = 0 ,
\end{equation}
for Love waves, and
\begin{equation}\label{Rayleigh1_ti}
-\frac{\partial}{\partial Z}C_{1313}\frac{\partial\varphi_1}{\partial Z}-\ii|\xi|\left(\frac{\partial}{\partial Z}(C_{1313}\varphi_3)+C_{1133}\frac{\partial}{\partial Z}\varphi_3\right)+C_{1111}|\xi|^2\varphi_1=\Lambda\varphi_1 ,
\end{equation}
\begin{equation}\label{Rayleigh2_ti}
-\frac{\partial}{\partial Z}C_{3333}\frac{\partial\varphi_3}{\partial Z}-\ii|\xi|\left(\frac{\partial}{\partial Z}(C_{1133}\varphi_1)+C_{1313}\frac{\partial}{\partial Z}\varphi_1\right)+C_{1313}|\xi|^2\varphi_3=\Lambda\varphi_3,
\end{equation}

\begin{eqnarray}
\ii|\xi| \vphi_3(0) + \pdpd{\vphi_1}{Z}(0) &=& 0 ,\label{Rayleighboundary1_ti}
\\
\ii C_{1133} |\xi| \vphi_1(0) + C_{3333} \pdpd{\vphi_3}{Z}(0)\label{Rayleighboundary2_ti}
&=& 0 ,
\end{eqnarray}
for Rayleigh waves.

\subsection{Directional decoupling}

We assume that $\xi = (|\xi|,0,0)$, that is, the phase direction of
propagation is $x_1$, and assume that the medium is monoclinic with
symmetry plane $(x_1,Z)$. Then $C_{ijkl}$ is determined by 13 nonzero
components:
\[
   C_{1111}, ~C_{1122}, ~C_{1133}, ~C_{2222}, ~C_{2233}, ~C_{3333},
   C_{1223}, ~C_{2323}, ~C_{1212}, ~C_{1113}, ~C_{1333}, ~C_{1313},
   ~C_{2213} .
\]
The symmetry plane, $(x_1,Z)$, is called the sagittal plane. Then
\begin{multline*}
H_{0}(x,\xi)u=\left[ \left(\begin{array}{ccc}
-\frac{\partial}{\partial Z}(C_{1313}\frac{\partial}{\partial Z})& 0 & -\frac{\partial}{\partial Z}(C_{1333}\frac{\partial}{\partial Z})\\
0 & -\frac{\partial}{\partial Z}(C_{2323}\frac{\partial}{\partial Z}) & 0\\
-\frac{\partial}{\partial Z}(C_{1333}\frac{\partial}{\partial Z}) & 0 & -\frac{\partial}{\partial Z}(C_{3333}\frac{\partial}{\partial Z})
\end{array}\right) \right.\\
-\ii|\xi| \left(\begin{array}{ccc}
C_{1113}\frac{\partial}{\partial Z} & 0 & C_{1133}\frac{\partial}{\partial Z}\\
0 & C_{1223}\frac{\partial}{\partial Z} & 0\\
C_{1313}\frac{\partial}{\partial Z} & 0 & C_{1333}\frac{\partial}{\partial Z}
\end{array}\right)
-\ii|\xi| \left(\begin{array}{ccc}
\frac{\partial}{\partial Z}(C_{1113}) & 0 & \frac{\partial}{\partial Z}(C_{1313})\\
0 & \frac{\partial}{\partial Z}(C_{1223}) & 0\\
\frac{\partial}{\partial Z}(C_{1133}) & 0 & \frac{\partial}{\partial Z}(C_{1333})
\end{array}\right)\\
\left.
-\ii|\xi| \left(\begin{array}{ccc}
C_{1113}\frac{\partial}{\partial Z}& 0 &C_{1313} \frac{\partial}{\partial Z}\\
0 &C_{1223} \frac{\partial}{\partial Z} & 0\\
C_{1133}\frac{\partial}{\partial Z} & 0 & C_{1333}\frac{\partial}{\partial Z}
\end{array}\right)+|\xi|^2\left(\begin{array}{ccc}
C_{1111} & 0 & C_{1113}\\
0 & C_{1212} & 0\\
C_{1113} & 0 & C_{1313}
\end{array}\right)\right]\left(\begin{array}{c}u_1\\ u_2\\u_3\end{array}\right) .
\end{multline*}
We obtain the equations for the eigenfunctions:
\begin{eqnarray}\label{monoclinic_first}
   -\frac{\partial}{\partial Z}C_{2323}\frac{\partial}{\partial Z}u_2-\ii C_{1223}|\xi|\frac{\partial}{\partial Z}u_2-\ii |\xi|\frac{\partial}{\partial Z}(C_{1223}u_2)+|\xi|^2C_{1212}u_2 &=& \Lambda u_2 ,
\\
C_{2323}\frac{\partial u_2}{\partial Z}(0)-\ii |\xi|C_{1223}u_2(0) &=&0\label{monoclinic_love_final}
\end{eqnarray}
for Love waves, and
\begin{multline}
-\frac{\partial}{\partial Z}C_{1313}\frac{\partial}{\partial Z}u_1-\frac{\partial}{\partial Z}C_{1333}\frac{\partial}{\partial Z}u_3-\ii|\xi|C_{1113}\frac{\partial}{\partial Z}u_1-\ii|\xi|C_{1133}\frac{\partial}{\partial Z}u_3
\\
\ii|\xi|\frac{\partial}{\partial Z}(C_{1113}u_1)-\ii|\xi|\frac{\partial}{\partial Z}(C_{1313}u_3)+|\xi|^2C_{1111}u_1+|\xi|^2C_{1113}u_3=\Lambda u_1
\end{multline}
\begin{multline}
-\frac{\partial}{\partial Z}C_{1333}\frac{\partial}{\partial Z}u_1-\frac{\partial}{\partial Z}C_{3333}\frac{\partial}{\partial Z}u_3-\ii|\xi|C_{1313}\frac{\partial}{\partial Z}u_1-\ii|\xi|C_{1333}\frac{\partial}{\partial Z}u_3
\\
\ii|\xi|\frac{\partial}{\partial Z}(C_{1133}u_1)-\ii|\xi|\frac{\partial}{\partial Z}(C_{1333}u_3)+|\xi|^2C_{1113}u_1+|\xi|^2C_{1313}u_3=\Lambda u_3
\end{multline}
\begin{equation}
C_{1313}\frac{\partial u_1}{\partial Z}(0)+\ii|\xi|C_{1113}u_1(0)+C_{1333}\frac{\partial u_3}{\partial Z}(0)+\ii|\xi|C_{1313}u_3(0)=0 ,
\end{equation}
\begin{equation}\label{monoclinic_final}
C_{1333}\frac{\partial u_1}{\partial Z}(0)+\ii|\xi|C_{1133}u_1(0)+C_{3333}\frac{\partial u_3}{\partial Z}(0)+\ii|\xi|C_{1333}u_3(0)=0 ,
\end{equation}
for Sloping-Rayleigh waves. We observe that the surface-wave modes
decouple into Love waves and Sloping-Rayleigh waves when the direction of propagation is in the saggital plane.

For ``global'' decoupling, decoupling takes place in every direction
of $\xi$. Thus transverse isotropy is the minimum symmetry to obtain
global decoupling; for a further discussion, see \cite{Anderson} . For
more generally anisotropic media, decoupling is no longer
possible. Then there only exist ``generalized'' surface-wave modes
propagating with elliptical particle motion in three dimensions
\cite{Crampin0, Crampin1}. Again, only in particular directions of
saggital symmetry, however, the Sloping-Rayleigh and Love waves
separate \cite{Crampin2} up to the leading order term. This loss of
polarization can be used to recognize anisotropy of the medium.

\section{Weyl's laws for surface waves}
\label{sec:6}

Weyl's law was first established for an eigenvalue problems for the
Laplace operator,
\[
\begin{cases}
   -\Delta  u &= \lambda u~~~\text{in}~~\Omega ,
\\
            u &= 0~~~\text{on}~~\partial\Omega,
\end{cases}
\]
for $\Omega$ a two- or three-dimensional bounded domain. One defines
the counting function, $N(b) := \#\{\lambda_n\leq b\}$, where
$\lambda_n$ are the eigenvalues arranged in a non-decreasing order
with counting multiplicity. Weyl \cite{Weyl1,Weyl2,Weyl3} proved the
following behaviors,
\begin{equation}\label{LapWeyl1}
   N(b) = \frac{|\Omega|}{4\pi}b+o(b) ,\quad
   b \rightarrow \infty ,
\end{equation}
for $\Omega \subset \mathbb{R}^2$, and
\begin{equation}\label{LapWeyl2}
   N(b) = \frac{|\Omega|}{6\pi^2}b^{3/2}+o(b^{3/2}) ,\quad
   b \rightarrow \infty ,
\end{equation}
for $\Omega \subset \mathbb{R}^3$. Formulae of these type are referred
to as Weyl's laws. For a rectangular domain, one can easily find
explicitly the eigenvalues and verify the formulae. These imply that
the leading order asymptotics of the number of eigenvalues is
determined by the area/volume of the domain only. Extensive work has
been done generalizing this result and deriving expressions for the
remainder term. We refer to \cite{ANPS} for a comprehensive
overview. Here, we present Weyl's laws for surface waves.

\subsection{Isotropic case}

\subsubsection*{Love waves}

 First, we consider the equations
(\ref{love}) and (\ref{love_b}) for Love waves in isotropic media. We
use the notation $C_S^2(x,Z) := \hat{\mu}(x,Z)$.  We have the
following Hamiltonian,
\begin{equation}\label{semiHelmholtz}
   H_0(x,\xi) v := -\frac{\partial}{\partial Z} \left(
         C_S^2\frac{\partial }{\partial Z}v \right)
                   + C_S^2 |\xi|^2 v ,
\end{equation}
where $H_0(x,\xi)$ is an ordinary differential operator in $Z$, with domain
\[
   \mathcal{D} = \left\{ v \in H^2(\mathbb{R}^-)\ :\
                \frac{\mathrm{d}v}{\mathrm{d}Z}(0)=0\right\} .
\]
We make use of results in \cite{CDV2}. Although in \cite{CDV2}, the
boundary condition is of Dirichlet type, the spectral property is
essentially the same. A straightforward extension of the argument in
\cite{Gladwell} suffices:
\begin{gather*}
   \sigma(H_0(x,\xi)) = \sigma_{p}(H_0(x,\xi))
                            \cup \sigma_{c}(H_0(x,\xi)) , 
\end{gather*}
where the point spectrum $\sigma_{p}(H_0(x,\xi))$, consists of a
finite set of eigenvalues in
\[
   \left( \inf_{Z \leq 0} C_S(x,Z)^2 |\xi|^2, C_S(x,Z_I)^2 |\xi|^2
          \right)
\]
and the continuous spectrum is given by
\[
   \sigma_{c}(H_0(x,\xi)) = [C_S(x,Z_I)^2 |\xi|^2,\infty) .
\]
Since the problem at hand is scalar and one-dimensional, the
eigenvalues are simple. We note that, here, the essential spectrum is
purely the continuous spectrum.

\begin{figure}[h]
\centering
\includegraphics[width=2 in]{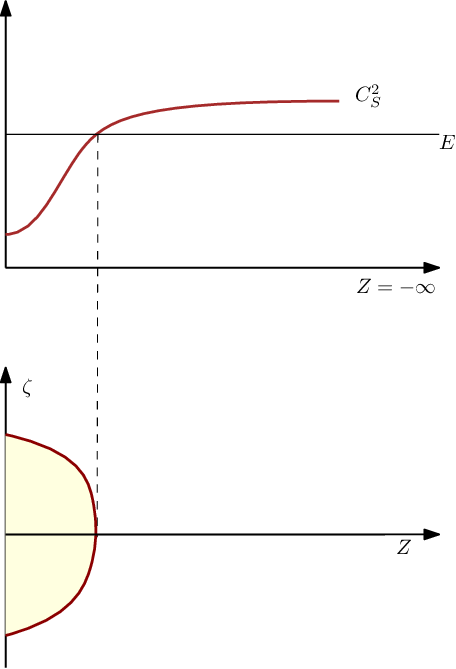}
\caption{Illustration of Weyl's law for Love waves. The bottom figure
  indicates the relevant volume.}
\label{fig:Weyl}
\end{figure}

We write $N(x,\xi,E) = \#\{\Lambda_\alpha(x,\xi)\leq E|\xi|^2\}$,
where $\Lambda_\alpha(x,\xi)$ is an eigenvalue of $H_0(x,\xi)$. Weyl's
law gives a quantitative asymptotic approximation of $N(x,\xi,E)$ in
terms of $|\xi|$ (Figure \ref{fig:Weyl}).

\medskip\medskip

\begin{theorem}[Weyl's law for Love waves]
For any $E < C^2_S(x,Z_I)$, we have 
\begin{equation}
   N(x,\xi,E) = \frac{|\xi|}{2\pi}
        \left( |\{(Z,\zeta):\ C^2_S(x,Z) (1+\zeta^2) \leq E \}|
                + o(1) \right),
\end{equation}
where $|\cdot|$ denotes the measure of the set.
\end{theorem}

\medskip\medskip

The classical Dirichlet-Neumann bracketing technique could be used to
prove this theorem \cite{RS}. However, we can also adapt the proof of
Theorem~\ref{WeylRayleigh} below.

We observe that $\Lambda_\alpha(x,\xi)$ is a smooth function of
$(x,\xi) \in T^*\mathbb{R}^2$ for all $1 \leq \alpha \leq
N(x,\xi,E)$. Moreover, the corresponding eigenfunctions,
$\varphi_\alpha(x,\xi,Z)$, decay exponentially in $Z$ for large
$Z$. More precisely, let $Z_b = \sup\{Z\ |\ C^2_S(x,Z) = b\}$, then
for all $\Lambda_\alpha < b$, the corresponding $\varphi_\alpha$ are
exponentially decaying from $Z = Z_b$ to $-\infty$. Hence, the energy
of $\varphi_\alpha$ is concentrated near a neighborhood of
$[Z_b,0]$. As $b \rightarrow Z_0(x)$, then the energy of
$\varphi_\alpha$ with $\Lambda_\alpha < b$ becomes more and more
concentrated near the boundary. Thus those low-lying eigenvalues
correspond to surface waves.

\medskip\medskip

\begin{remark}
For Love waves in a monoclinic medium described by equations
$(\ref{monoclinic_first})-(\ref{monoclinic_love_final})$, Weyl's law
reads: For any $E <C_{1212}-\frac{4C_{1223}^2}{C_{2323}} $ , we have
\begin{equation}
   N(x,\xi,E) = \frac{|\xi|}{2\pi}
        \left( |\{(Z,\zeta):
     C_{2323}\zeta^2+2C_{1223}\zeta+C_{1212}\leq E\}| + o(1) \right) .
\end{equation}
\end{remark}

\subsubsection*{Rayleigh waves}

Here, we establish Weyl's law for Rayleigh waves,
cf.~(\ref{Rayleigh1})-(\ref{Rayleighboundary2}). Now,
\begin{equation}\label{Hamiltonian_Rayleigh}
H_0(x,\xi)\left(\begin{array}{c}
\varphi_1\\
\varphi_3
\end{array}\right)
=\left(\begin{array}{c}
-\frac{\partial}{\partial Z}(\hat{\mu}\frac{\partial\varphi_1}{\partial Z})-\ii|\xi|\left(\frac{\partial}{\partial Z}(\hat{\mu}\varphi_3)+\hat{\lambda}\frac{\partial}{\partial Z}\varphi_3\right)+(\hat{\lambda}+2\hat{\mu})|\xi|^2\varphi_1\\
-\frac{\partial}{\partial Z}\left((\hat{\lambda}+2\hat{\mu})\frac{\partial\varphi_3}{\partial Z}\right)-\ii|\xi|\left(\frac{\partial}{\partial Z}(\hat{\lambda}\varphi_1)+\hat{\mu}\frac{\partial}{\partial Z}\varphi_1\right)+\hat{\mu}|\xi|^2\varphi_3
\end{array}\right) .
\end{equation}
We use $\zeta$ to denote the Fourier variable for
$\frac{1}{\ii|\xi|}\frac{\partial}{\partial Z}$. For fixed $(x,\xi)$,
we view $H_0(x,\xi)$ as a semiclassical pseudodifferential operator
with $\frac{1}{|\xi|}$ as semiclassical parameter, $h$, as before. We
denote $h_0(x,\xi)=h^2H_0(x,\xi)$. Then the principal symbol for
$h_0(x,\xi)$ is
\[h_0(x,\xi)(Z,\zeta)=\left[\begin{array}{cc}  \hat{\mu}\zeta^2& 0\\
0 &(\hat{\lambda}+2\hat{\mu})\zeta^2\end{array}\right]+\left[\begin{array}{cc}  0 &(\hat{\lambda}+\hat{\mu})\zeta\\
(\hat{\lambda}+\hat{\mu})\zeta& 0\end{array}\right]+\left[\begin{array}{cc} \hat{\lambda}+2\hat{\mu}& 0\\
0 &\hat{\mu}\end{array}\right] .
\]
There exists $\Pi(\zeta)$, namely, 
\[
   \Pi(\zeta) = \left[\begin{array}{cc} 1 & -\zeta \\ \zeta & 1
                    \end{array}\right] ,
\]
such that
\[
   h_0(x,\xi)(Z,\zeta) \, \Pi = \Pi
   \left[\begin{array}{cc} (\hat{\lambda}+2\hat{\mu})(1+\zeta^2) & 0\\
         0 & \hat{\mu}(1+\zeta^2) \end{array}\right] .
\]
We write $C_P^2(x,Z) = (\hat{\lambda}+2\hat{\mu})(x,Z)$. We obtain

\medskip\medskip

\begin{theorem}[Weyl's law for Rayleigh waves]\label{WeylRayleigh}
For any $E<\hat{\mu}(x,Z_I)$, define
\[
   N(x,\xi,E) = \#\{ \Lambda_\alpha:
                              \Lambda_\alpha \leq E |\xi|^2 \} ,
\]
where $\Lambda_\alpha$ are the eigenvalues of $H_0(x,\xi)$. Then

\begin{equation}\label{sharp Weyl}
   N(x,\xi,E) = \frac{|\xi|}{2\pi }
        \left(\left|\left\{ (Z,\zeta): C_S^2(x,Z)(1+\zeta^2)
                   \leq E \right\}\right|
    + \left|\left\{  (Z,\zeta):
               C_P^2(x,Z)(1+\zeta^2) \leq E \right\}\right|
                  + o(1) \right) .
\end{equation}
\end{theorem}

\medskip\medskip

\begin{proof}
Let $\chi_1(Z),\ \chi_2(Z) \in C^{\infty}(\mathbb{R}^-)$ be
real-valued, non-negative functions, such that $\chi_1(Z) = 0$ for $Z
\geq -\frac{\delta}{2}$, $\chi_1(Z) = 1$ for $Z \leq -\delta$, and
$\chi^2_1(Z) + \chi^2_2(Z) = 1$ for some $\delta > 0$. For any
non-negative $f \in C_c^\infty(\mathbb{R})$, $\operatorname{supp} f
\subset (-\infty, \hat{\mu}(x,Z_I))$, $f(h^2 H_0)$ is a trace class
operator; we have
\[
   \operatorname{trace} f(h^2 H_0)
        = \operatorname{trace}\left(\chi_1^2 f(h^2 H_0)\right)
             + \operatorname{trace}\left(\chi_2^2 f(h^2 H_0)\right) .
\]
We analyze the two terms successively. We note that
\[
   \operatorname{trace}\left(\chi_1^2 f(h^2 H_0)\right)
        =\operatorname{trace}\left(\chi_1 f(h^2 H_0) \chi_1 \right) ,
\]
and that $\chi_1 f(h^2 H_0)\chi_1$ is a pseudodifferential operator on
$\mathbb{R}$ with principal symbol,
\[
   \chi_1(Z) f(h_0 (Z,\zeta)) \chi_1(Z') .
\]
Revisiting the diagonalization of $h_0$ above, let $\mathfrak{U}$  be the pseudodifferential
operator with symbol
\[
   \frac{1}{\sqrt{1+\zeta^2}} \Pi(\zeta) ;
\]
it follows immediately that $ \mathfrak{U}$ is unitary on $L^2(\mathbb{R})$. Then
\[
   \chi_1f(h^2 H_0) \chi_1 = \chi_1 \mathfrak{U} f(D_0) \mathfrak{U}^* \chi_1 + O(h)
               = \mathfrak{U} \chi_1 f(D_0) \chi_1 \mathfrak{U}^* + O(h) ,
\]
with $D_0$ being the pseudodifferential operator with symbol
\[
   \left[\begin{array}{cc}
      (\hat{\lambda}+2\hat{\mu}) (1+\zeta^2) & 0 \\
               0 & \hat{\mu}(1+\zeta^2) \end{array}\right] .
\]

Let $D_P$ and $D_S$ be the pseudodifferential operators with symbols
$(\hat{\lambda}+2\hat{\mu}) (1+\zeta^2)$ and $\hat{\mu} (1+\zeta^2)$,
respectively. Based on the discussion at the end of Appendix \ref{App_Opth},
\begin{multline}\label{chi1part}
   \operatorname{trace}\left(\chi_1 f (h^2 H_0) \chi_1\right)
      = \operatorname{trace}(\chi_1 f(D_0) \chi_1) + O(h)
\\[0.25cm]
      = \operatorname{trace}(\chi_1 f(D_S) \chi_1) +
              \operatorname{trace}(\chi_1 f(D_P) \chi_1) + O(h)
\\
      = \frac{1}{2\pi h} \left(\int_{-\infty}^\infty
             \int_{-\infty}^0 \chi_1(Z)^2
         \left[f(\hat{\mu}(Z) (1+\zeta^2))
              + f((\hat{\lambda}+2\hat{\mu})(Z) (1+\zeta^2)) \right]
           \mathrm{d}Z \mathrm{d}\zeta + O(h)\right) .
\end{multline}
We then observe that
\[
   0 \leq \operatorname{trace} \left( \chi_2^2 f(h^2 H_0) \right)
              \leq \operatorname{trace}f(h^2 H_{0,\delta}) ,
\]
where $H_{0,\delta}$ is the operator $H_0$ restricted to
$[-\delta,0]$. To show this, we choose the eigenvalues, $\mu_j$, and
associated unit eigenvectors, $\mathbf{e}_j$, of $h^2H_{0,\delta}$ subject
to the Neumann boundary condition. Then
$\{\mathbf{e}_j\}_{j=1}^\infty$ forms an orthonormal basis of
$L^2([-\delta,0])$. We also consider an arbitrary orthonormal basis
$\{\mathbf{f}_j\}_{j=1}^\infty$ for $L^2([-\infty,-\delta])$. Clearly,
$\{\mathbf{e}_i,\mathbf{f}_j\}_{i,j=1}^\infty$ form an orthonormal
basis for $L^2(\mathbb{R}^-)$. Therefore,
\begin{multline}
   \operatorname{trace}\left( \chi_2^2 f(h^2 H_0) \right)
      = \sum_{j=1}^\infty \langle\chi_2^2 f(h^2 H_0) \mathbf{e}_j,
              \mathbf{e}_j\rangle
      + \sum_{j=1}^\infty \langle\chi_2^2 f(h^2 H_0) \mathbf{f}_j,
              \mathbf{f}_j\rangle
\\
   = \sum_{j=1}^\infty \langle\chi_2^2 f(h^2 H_0) \mathbf{e}_j,
             \mathbf{e}_j\rangle
   = \sum_{j=1}^\infty \langle\chi_2^2 f(\mu_j) \mathbf{e}_j,
             \mathbf{e}_j\rangle
     \leq \sum_{j=1}^\infty \langle f(\mu_j) \mathbf{e}_j,
             \mathbf{e}_j\rangle
   = \operatorname{trace} f(h^2 H_{0,\delta}) .
\end{multline}

For every small
$\epsilon>0$, we can take $f(u) \leq
\mathbf{1}_{[-\epsilon,E+\epsilon]}(u)$ , where $\mathbf{1}_{[-\epsilon,E+\epsilon]}(u)$ is the indicator function of $[-\epsilon,E+\epsilon]$, that is,
\[
\mathbf{1}_{[-\epsilon,E+\epsilon]}(u)=
\begin{cases}
1,~~~~\text{ if }~u\in [-\epsilon,E+\epsilon]\\
0,~~~~\text{ elsewhere}.
\end{cases}\]
Then
\begin{equation}\label{chi2part}
\begin{split}
   \operatorname{trace} f(h^2 H_{0,\delta})
       &\leq \#\{E(h)\ |\ E(h) \leq E +\epsilon,\
              E(h)\text{ is an eigenvalue of }h^2 H_{0,\delta}\}
\\
       &\leq \frac{C \delta}{h} .
\end{split}\end{equation}
To show this, we rescale $z = \frac{Z}{h}$. Then for $\epsilon<1$,
\begin{multline}
   \#\{E(h): E(h) \leq E+\epsilon,
            E(h)\text{ is an eigenvalue of }h^2 H_{0,\delta}\}
\\
   \leq \#\{\lambda\ |\ \lambda \leq E+1,
           \lambda\text{ is an eigenvalue of }A\},
\end{multline}
where
\begin{equation}\label{A}
   A = \left(\begin{array}{cc}
   -\frac{\partial}{\partial z} (\hat{\mu}
        \frac{\partial}{\partial z}\cdot)+ (\hat{\lambda}+2\hat{\mu})\cdot&
  -\ii \left(\frac{\partial}{\partial z} (\hat{\mu}\cdot)
           + \hat{\lambda}\frac{\partial}{\partial z}\cdot\right)
\\
    -\ii (\frac{\partial}{\partial z} (\hat{\lambda}\cdot)
           + \hat{\mu}\frac{\partial}{\partial z}\cdot) &
   -\frac{\partial}{\partial z}\left( (\hat{\lambda}+2\hat{\mu})
           \frac{\partial}{\partial z}\cdot\right) +\hat{\mu}\cdot 
   \end{array}\right)
\end{equation}
on $[-\frac{\delta}{h},0]$ with the Neumann boundary condition
applied. Without loss of generality, we can assume that
$\frac{\delta}{h}$ is an integer. Then we divide the interval
$[-\frac{\delta}{h},0]$ into $\frac{\delta}{h}$ intervals of the same
length, $1$, and let $A_j$ be operator $A$ restricted to $[-j,-j+1]$,
with the Neumann boundary condition applied. Then, as in the proof of
Theorem \ref{existence},
\[
   A \geq \bigoplus_{i=1}^{\frac{\delta}{h}} A_j .
\]
The number of eigenvalues of each $A_j$ below $E+1$ can be bounded by
a constant since the corresponding quadratic forms have a uniform
lower bound. Hence, we obtain (\ref{chi2part}).

Finally, we combine (\ref{chi1part}) with (\ref{chi2part}), let
$\delta \rightarrow 0$, and note that $\chi_1 \rightarrow
\mathbf{1}_{[0,\infty)}$, so that
\begin{equation}\label{traceformula}
   \operatorname{trace}\left(f(h^2H_0)\right) = \frac{1}{2\pi h}
       \left(\int_{-\infty}^\infty \int_{-\infty}^0
   \left[f(\hat{\mu}(Z) (1+\zeta^2))
       + f((\hat{\lambda}+2\hat{\mu})(Z) (1+\zeta^2))\right]
              \mathrm{d}Z\mathrm{d}\zeta+\mathcal{O}(h)\right) .
\end{equation}

Using a technique similar to the one used in the proof of
\cite[Theorem 14.11]{Zworski}, we construct $f_1^\epsilon,
f_2^\epsilon \in C_c^\infty(\mathbb{R})$ by regularizing
$\mathbf{1}_{[0,E]}$, with $f_1^\epsilon, f_2^\epsilon \leq
\mathbf{1}_{[-\epsilon,E+\epsilon]}$, such that
\[
   f_1^\epsilon(u) \leq \mathbf{1}_{[0,E]}(u) \leq f_2^\epsilon(u) .
\]
More presicely, we may take
\[J(u)=\begin{cases}
k\exp[-1/(1-|u|^2)]~~~~~&\text{if}~|u|<1\\
0~~~~~~~~~~~~~~~~~~~~~~~~~~~~~~&\text{if}~|u|\geq 1,
\end{cases}\]
where $k$ is chosen so that $\int_{\mathbb{R}}J(u)\mathrm{d}u=1$. Let
$J_\epsilon(u)=\frac{1}{\epsilon}J(\frac{u}{\epsilon})$, and
\begin{eqnarray*}
   f_1^\epsilon &=& J_\epsilon * \mathbf{1}_{[\epsilon,E-\epsilon]} ,
\\
   f_2^\epsilon &=& J_\epsilon*\mathbf{1}_{[-\epsilon,E+\epsilon]} .
\end{eqnarray*}
While observing that $f_1^\epsilon, f_2^\epsilon \rightarrow
\mathbf{1}_{[0,E]}$, and using the estimates above, taking $\epsilon
\rightarrow 0$ completes the proof.
\end{proof}

\subsection{Anisotropic case}

We extend Weyl's law from isotropic to anisotropic media. The proof
for Rayleigh waves can be naturally adapted. We write the eigenvalues
of the symmetric matrix-valued symbol defined in (\ref{principle}) as
$\mathfrak{C}_1(x,\hat{\xi},Z,\zeta)$,
$\mathfrak{C}_2(x,\hat{\xi},Z,\zeta)$ and
$\mathfrak{C}_3(x,\hat{\xi},Z,\zeta)$

\medskip\medskip

\begin{theorem}
Assume the three eigenvalues
$\mathfrak{C}_i(x,\hat{\xi},Z,\zeta),i=1,2,3$ are smooth in
$(Z,\zeta)$, then for any $E < v_L(x,\hat{\xi},Z_I)|\xi|^2$, let
\[
   N(x,\xi,E) = \#\{\Lambda_\alpha\ :\
                   \Lambda_\alpha(x,\xi) \leq E |\xi|^2\} .
\]
Then
\[
   N(x,\xi,E) = \frac{|\xi|}{2\pi}
       (|\{(Z,\zeta)\ :\ \mathfrak{C}_1 \leq E\}|
                  +|\{(Z,\zeta)\ :\ \mathfrak{C}_2 \leq E\}|
                        +|\{(Z,\zeta)\ :\ \mathfrak{C}_3 \leq E\}| + o(1)) .
\]
\end{theorem}

\section{Surface waves as normal modes}
\label{sec:7}

\begin{figure}[h]
\centering
\includegraphics[width=4 in]{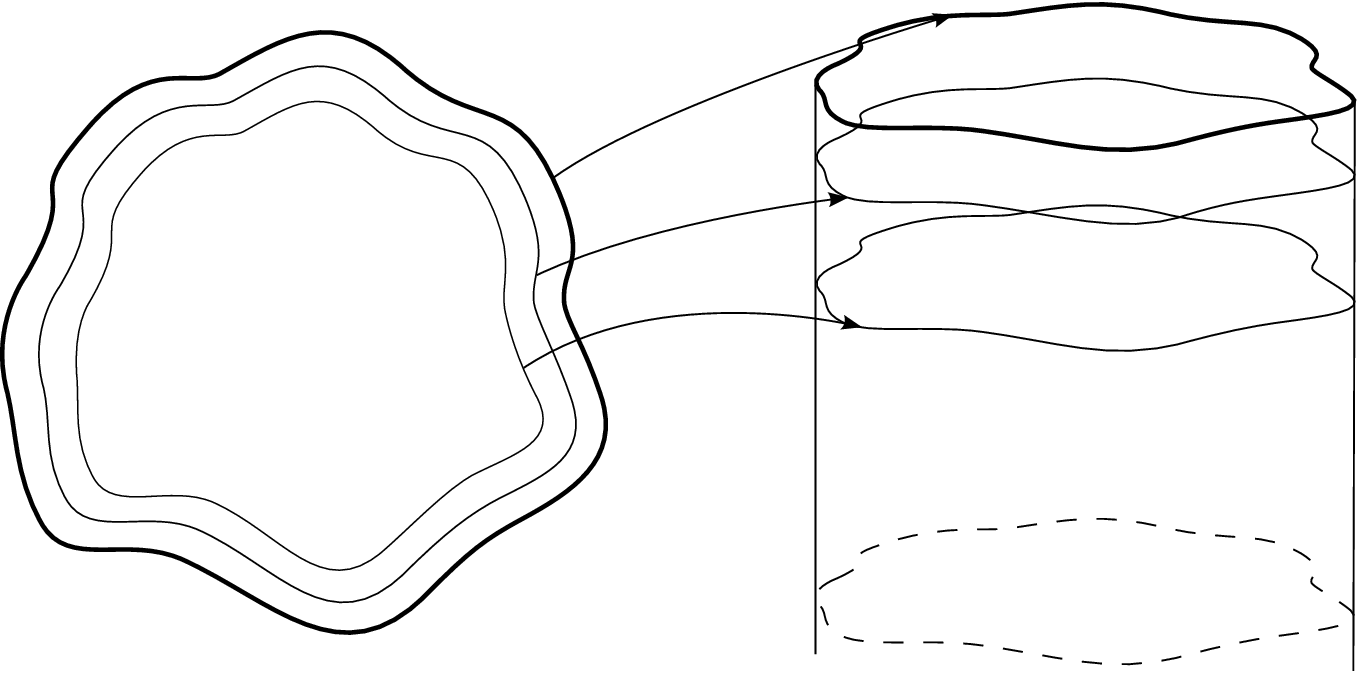}
\caption{Illustration of the diffeomorphism $\phi$.}
\label{product}
\end{figure}

In this section, we identify surface wave modes with normal modes.  We consider
the Earth as a unit ball $B_1$. There is a global diffeomophism, $\phi$,
with
\begin{equation*}
\begin{split}
   \phi :\ B_1 \setminus \{0\}&\rightarrow S^2 \times \mathbb{R}^- ,
\\
    \phi(B_r)&= S^2 \times \left\{1-\frac{1}{r}\right\} ,~r \neq 0 ,
\end{split}
\end{equation*}
where $S^2$ is the unit sphere centered at the origin.
For an open and bounded  subset $U \subset S^2$, the cone
region, $\{(\Theta,r)\ |\ \Theta \in U ,\ 0<r<1\}$, is diffeomorphic
to $U \times \mathbb{R}^-$; hence, we can find global coordinates for
$U$ and we may consider our system on the domain $S^2 \times
\mathbb{R}^-$. More generally, we consider the system on any
Riemannian manifold of the form $M = \partial M \times \mathbb{R}^-$
with metric
\[
   g = \left(\begin{array}{cc}g' & 0 \\0 & 1\end{array}\right) .
\]
Thus $\partial M$ is a compact Riemannian manifold with metric
$g'$. For a ``nice'' domain $\Omega$, a neighborhood of the boundary
is diffeomorphic to $M$, where the metric $g'$ is the induced metric
of the boundary of $\Omega$. We note that $g'$ captures possible
``topography'' of $\partial\Omega$ (Figure \ref{product}).

Let $x'=(x_1,x_2)$ be local coordinates on $\partial M$, and $x =
(x_1,x_2,x_3)=(x',Z)$ be local coordinates on $M$. Then $g$ can be
represented by
\[g=\left(
\begin{array}{ccc}
g_{11}&  g_{12} &0\\
g_{21}& g_{22}&0\\
0& 0 &1
\end{array}\right) .
\] 
We represent $\hat{H}$ on $M$ in local coordinates, similarly as in
$(\ref{Full})$.  Now the displacement $u_l$ is a covariant vector,
\begin{multline}\label{Full_tensorial}
\hat{H}_{i}^{~l} u_l=
-\frac{\partial}{\partial Z} C_i^{~33l}(x,Z)
\frac{\partial}{\partial Z}u_l
\\
- \epsilon \sum_{j=1}^{2}
C_i^{~j3l}(x,Z)\nabla_j
\frac{\partial}{\partial Z}u_l
- \epsilon \sum_{k=1}^{2}
C_i^{~3kl}(x,Z) \frac{\partial}{\partial Z}
\nabla_ku_l
- \epsilon \sum_{k=1}^{2}
\left( \frac{\partial}{\partial Z} C_i^{~3kl}(x,Z) \right)
\nabla_ku_l\\
- \epsilon^2 \sum_{j,k=1}^{2}
C_i^{~jkl}(x,Z)\nabla_j
\nabla_ku_l-\epsilon \sum_{j=1}^{2}
\left(\nabla_jC_i^{~j3l}(x,Z) \right)
\frac{\partial}{\partial Z}u_l\\
- \epsilon^2 \sum_{j,k=1}^{2}
\left(\nabla_jC_i^{~jkl}(x,Z) \right)
\nabla_ku_l.
\end{multline}
Here, $\nabla_k$ is the covariant derivative associated with $x_k$.
Normal modes can be viewed as solutions to the eigenvalue problem for
\begin{equation}\label{normal}
   \hat{H}_{i}^{~l} u_l =\omega^2 u_i,
\end{equation}
with Neumann boundary condition. In fact, (\ref{normal}) is
asymptotically equivalent to 
\begin{equation}\label{surfacemodes}
   a_{\alpha,\epsilon}(\cdot,\epsilon \nabla) \Psi(x)
                  = \omega^2 \Psi(x) ,
\end{equation}
where $a_{\alpha,\epsilon}$ is the pseudodifferential operator on
$\partial M$ defined in Theorem~\ref{conjugation_sys}.  By the
estimates, $\Lambda_\alpha(x,\xi) \geq C|\xi|^2$, it follows that
$a_{\alpha,\epsilon}(x,\xi)$ is elliptic implying that
\[
   \| a_{\alpha,\epsilon}(\cdot,\epsilon \nabla) \Psi
            \|_{L^2(\partial M)} \geq C \|\Psi\|_{H^2(\partial M)} .
\]
Thus $a_{\alpha,\epsilon}(\cdot,\epsilon \nabla)^{-1}$ is compact, and
hence the spectrum of $a_{\alpha,\epsilon}(\cdot,\epsilon \nabla)$ is
discrete.

Now, let $\omega$ be an eigenfrequency, then we construct asymptotic
solutions of (\ref{surfacemodes}) of the form
\begin{equation}
   \Psi = \left( \sum_{k=0}^\infty \epsilon^k \mathcal{B}_k(x)
             \right) e^{\ii \frac{\psi}{\epsilon}} .
\end{equation}
Inserting this expression into (\ref{surfacemodes}), from the
expansion, we find that
\[
   \Lambda_\alpha(x,\partial\psi) = \omega^2
\]
and
\begin{equation}
   L \mathcal{B}_0 = 0
\end{equation}
\begin{equation}
   L \mathcal{B}_k = F(\psi,\mathcal{B}_0,\cdots,\mathcal{B}_{k-1}) .
\end{equation}
Here, $L$ and $F$ are defined similar to those defined in Section
\ref{param}. Thus
\[
   u = \operatorname{Op}_\epsilon(\Phi_{\alpha,\epsilon}) (\Psi)
\]
will be an asymptotic solution of (\ref{normal}).

\medskip\medskip

\begin{remark}[Radial manifolds]
Under the assumption of transverse isotropy and lateral homogeneity,
that is, fixing an $x$, and using the normal coordinates at $x$, we
have (\ref{love_ti})-(\ref{Rayleighboundary2_ti}) with
$C_{ijkl}(x,Z)=C_{ijkl}(Z)$. Then we can construct asymptotic
modes. Now $|\xi|^2$ corresponds to
$-\Delta_{g'}^\epsilon=-\epsilon^2\Delta_{g'}$, where $\Delta_{g'}$ is
the Laplacian on $\partial M$. Then we consider the eigenvalue
problem:
\begin{equation}\label{love_ti_ball}
-\pdpd{}{Z} C_{1313} \pdpd{}{Z}u_2(x,Z) \,
-  C_{1212}\Delta_{g'}^\epsilon u_2(x,Z) = \omega^2 u_2(x,Z)
\end{equation}
and
\begin{equation}\label{Rayleigh1_ti_ball}
\begin{split}
-\frac{\partial}{\partial Z}C_{1313}\frac{\partial}{\partial Z}u_1(x,Z)-\ii\left(\frac{\partial}{\partial Z}C_{1313}+C_{1133}\frac{\partial}{\partial Z}\right)\sqrt{-\Delta_{g'}^\epsilon}u_3(x,Z)\\
-C_{1111}\Delta_{g'}^\epsilon u_1(x,Z)=\omega^2 u_1(x,Z)
\end{split}
\end{equation}
\begin{equation}\label{Rayleigh2_ti_ball}
\begin{split}
-\frac{\partial}{\partial Z}C_{3333}\frac{\partial}{\partial Z}u_3(x,Z)-\ii\left(\frac{\partial}{\partial Z}C_{1133}+C_{1313}\frac{\partial}{\partial Z}\right)\sqrt{-\Delta_{g'}^\epsilon}u_1(x,Z)\\
-C_{1313}\Delta_{g'}^\epsilon u_3(x,Z)=\omega^2 u_3(x,Z)
\end{split}
\end{equation}
Then, if we have the eigenvalues and eigenfunctions of $-\Delta_{g'}^\epsilon$,
\[-\Delta_{g'}^\epsilon\Theta_n^l(x)=k^2\Theta_n^l(x),\]
and solutions $\varphi(Z,k)$ for the system (\ref{love_ti})-(\ref{Rayleighboundary2_ti}) with $|\xi|=k$, we find
\[u(k,x,Z)=\varphi(Z,k)\Theta_n^l(x)\]
as the solutions to (\ref{love_ti_ball})-(\ref{Rayleigh2_ti_ball}). In
spherically symmetric models of the earth, $\partial M = S^2$ and
$\Theta_n^l(x)$ are the spherical harmonics.
\end{remark}

 \section*{Acknowledgements}
 M.V.d.H. gratefully acknowledges support from the Simons Foundation MATH + X program and National Science Foundation grant DMS-1559587. A.I.  is grateful to Department of Computational and Applied Mathematics, Rice University for the hospitality during his stay where the part of research was carried out. G.N. was partially supported by Grant-in-Aid for Scientific Research (15K21766 and 15H05740) of the Japan Society for the Promotion of Science. J.Z. thanks Brian Kennett for invaluable discussions.
\appendix

\renewcommand{\theequation}{\Alph{section}.\arabic{equation}}

\setcounter{equation}{0}
\section{Semiclassical pseudodifferential operators}\label{AA}
~\\

Here, we give a summary of the basic definition and properties of
semiclassical pseudodifferential operators which are used in the main
text. Let $A(\cdot,\cdot) :\ T^*\mathbb{R}^n\rightarrow
\mathbb{C}^{m\times m}$ be a symbol that is smooth in $(x,\xi)$. We
say that $A \in S(k)$, with $k\in\mathbb{Z}$, if \[
   \forall \alpha, \beta \in \mathbb{N}^n ,\quad
   |D_x^\alpha D_{\xi}^\beta A(x,\xi)|
                \leq C_{\alpha,\beta}\langle \xi\rangle^k ,
\]
with $\langle \xi \rangle = \sqrt{1+|\xi|^2}$.

Let $A_j \in S(k)$ for $j=0,1,\cdots$, and $\epsilon\in (0,\epsilon_0]$
for some small $\epsilon_0 > 0$. One says that $A \in S(k)$ is
asymptotic to $\sum_{j=0}^\infty \epsilon^jA_j$ and writes
\[
   A \sim \sum_{j=0}^\infty \epsilon^j A_j~~\text{in } S_k
\]
if for any $N=1,2,\cdots$ 
\[
   |\partial^\alpha(A-\sum_{j=0}^{N-1} \epsilon^j A_j)|
          \leq C_{\alpha,N} \epsilon ^N \langle\xi\rangle^k, ~~~\text{for any}~\alpha=0,1,2,\cdots.
\]
We refer to $A_0$ as the principal symbol of $A$. A semiclassical
psudodifferential operator associated with $A$ is defined as follows

\medskip\medskip

\begin{definition}\label{weyl}
Suppose that $A(x,\xi)$ is a symbol. We define the semiclassical
pseudodifferential operator,
\begin{gather}\label{weyl-quant}
   \operatorname{Op}_{\epsilon}(A) w(x) = A(.,\epsilon D) w(x)
     = \frac{1}{(2\pi\epsilon)^n} \iint A\left(x,\xi\right)
           e^{\ii \frac{\langle \xi,x - y \rangle}{\epsilon}}
                  w(y) \, \rmd y \rmd\xi,
\end{gather}
for any $w :\ \mathbb{R}^n \rightarrow \mathbb{C}^{m}$, which is
compactly supported.
\end{definition}

\medskip\medskip

We have the following mapping property: For any $u \in H^s$,
\[
   \| A(x,\epsilon D) u \|_{H^{s+k}}
                   \leq C_s \epsilon^k \|u\|_{H^s} ,
\]
for some constant $C_s > 0$. Here, $s$ is an arbitrary real number,
and $H^s$ denotes the $L^2$ Sobolev space in $\mathbb{R}^n$ with
exponent $s$.

\medskip\medskip

If $A(x,\xi)$ and $B(x,\xi)$ are two symbols, then
\[
   \operatorname{Op}_{\epsilon}(A)
      \operatorname{Op}_{\epsilon}(B)
                        = \operatorname{Op}_{\epsilon}(C)
\]
with
\[
   C(x,\xi) \sim
       \sum_{\alpha\geq 0} \frac{\ii^{|\alpha|}\epsilon^{|\alpha|}}{
            \alpha!} D_{\xi}^\alpha A(x,\xi) D^\alpha_x B(x,\xi) .
\]
We use the notation $\circ$ for the composition of two symbols: $C = A
\circ B$.

\medskip\medskip

\begin{definition}
We call a family of functions (distributions)
$u = \{ u_\epsilon \}_{0<\epsilon\leq\epsilon_0}$ admissible if there
exist constants $k$ and $N\geq 0$ such that
\[
   \| \chi u_\epsilon \|_{H^k} = \mathcal{O}(\epsilon^{-N}) ,
\]
for any $\chi \in C_c^\infty(\mathbb{R}^n)$. We denote by
$\mathcal{A}(\mathbb{R}^n)$ the space of such families.
\end{definition}

\medskip\medskip

\begin{definition}
The wavefront set of $u_\epsilon$, denoted by $WF(u_\epsilon)$, of an
admissible family $u=\{u_\epsilon\}_{0<\epsilon\leq\epsilon_0}$ is the
closed subset of $T^*\mathbb{R}^n$ which is defined as
\[
   (x_0,\xi_0) \notin WF(u_\epsilon)
\]
if and only if $\exists \chi\in C_c^\infty(\mathbb{R}^n)$,
$\chi(x_0)\neq 0$, such that
\[\mathcal{F}_\epsilon(\chi u_\epsilon)(\xi)=\mathcal{O}(\epsilon^\infty)\]
for $\xi$ close to $\xi_0$. Here $\mathcal{F}_\epsilon$ is the
semiclassical Fourier transform:
\[\mathcal{F}_\epsilon u_\epsilon(\xi)=\frac{1}{(2\pi\epsilon)^{n/2}}\int_{\mathbb{R}^n}e^{-\frac{\ii\langle x,\xi\rangle}{\epsilon}}u_\epsilon(x)\mathrm{d}x.\]
\end{definition}

\medskip\medskip

\begin{definition}\label{A4}
Let $U$ be an open set in $T^*\mathbb{R}^n$. The space of
microfunctions, $\mathcal{M}(U)$, in $U$ is the quotient
\[
   \mathcal{M}(U) = \mathcal{A}(\mathbb{R}^n) /
            \{u_\epsilon : WF(u_\epsilon)\cap U=\emptyset\} .
\]
\end{definition}

\medskip\medskip

Let $p(x,\xi)$ be a symbol such that
\[
   \nabla p \neq 0~~\text{on}~\{p=0\} .
\]
Let $(x,\xi)$ be the solution to the Hamilton system,
\[
   \frac{\partial x_k(y,\eta,t)}{\partial t}
       = \frac{\partial p(x,\xi)}{
                \partial \xi_k} ,\
   \frac{\partial \xi_k(y,\eta,t)}{\partial t}
       =-\frac{\partial p(x,\xi)}{
                \partial x_k} ,\quad
   (x,\xi)|_{t=0} = (y,\eta) .
\]
We denote by $\mathrm{exp}(t H_p)$ the map such that
$\mathrm{exp}(tH_p)(y,\eta) = (x,\xi)$. The map $\mathrm{exp}(tH_p)$
is called the Hamiltonian flow of $p$. For each
$u_\epsilon \in \mathcal{A}(\mathbb{R}^n)$ solving
\[
   p(x,\epsilon D) u_\epsilon = f_\epsilon ,
\]
where $\{f_\epsilon\}_{0<\epsilon\leq\epsilon_0} \subset
L^2(\mathbb{R}^n)$, we have $WF(u_\epsilon) \setminus WF(f_\epsilon)$
is invariant under the Hamiltonian flow of $p$.

\section{Operator theory}\label{App_Opth}
~\\

\begin{definition}
Let $\mathcal{H}$ be a Hilbert space, $D$ a dense subspace of
$\mathcal{H}$, and $A$ be an unbounded linear operator defined on $D$.
\begin{enumerate}
\item The adjoint $A^*$ of $A$ is the operator whose domain is $D^*$, where
\[D^*=\{v\in H: |\langle Au,v\rangle|\leq C(v)\|u\|\text{ for all }u\in D\},\]
and
\[\langle A^*v,u\rangle=\langle v,Au\rangle.\]
\item $A$ is called self-adjoint if $D^*=D$ and $A^*=A$.
\item $A$ is called symmetric if $D\subset D^*$ and $Au = A^*u$ for all $u\in D$.
\end{enumerate}
\end{definition}
~\\

\begin{definition}
\begin{enumerate}
\item Let $A$ be an unbounded self-adjoint operator densely defined on $\mathcal{H}$, with domain $D$. The spectrum of $A$ is
\[\sigma(A)=\mathbb{R}\setminus\{\lambda\in\mathbb{R}: (A-\lambda)^{-1}:\mathcal{H}\rightarrow \mathcal{H}\text{ is bounded}\}.\]
\item The set of all $\lambda$ for which $A-\lambda$ is injective and has dense range, but is not surjective, is called the continuous spectrum of $A$, denoted by $\sigma_{c}(A)$.
\item If there exists a $u\in D$ satisfying that $Au=\lambda u$, $\lambda$ is called an eigenvalue of $A$. The set of all eigenvalues is called the point spectrum of $A$, denoted by $\sigma_{p}(A)$.
\item The discrete spectrum $\sigma_{disc}(A)$ of $A$ is the set of eigenvalues of $A$ that have finite dimensional eigenspaces.
\item The essential spectrum $\sigma_{ess}(A)$ of $A$ is $\sigma(A)\setminus \sigma_{disc}(A)$.
\end{enumerate}
\end{definition}
~\\

\begin{theorem}\label{eigcri}
A number $\Lambda$ is in the essential spectrum of $A$ if and only if
there exists a sequence $\{u_k\}$ in $\mathcal{H}$ (called singular
sequence or Weyl sequence) such that
\begin{enumerate}
\item $\lim_{k\rightarrow\infty}\langle u_k,v\rangle=0,$ for all $v\in H;$

\item $\|u_k\|=1;$
\item $u_k\in D;$ 
\item $\lim_{k\rightarrow\infty}\|(A-\Lambda) u_k\|=0.$
\end{enumerate}
\end{theorem}
~\\

\textbf{Quadratic form}. A quadratic form is a map $q:Q(q)\times
Q(q)\rightarrow \mathbb{C}$, where $Q(q)$ is a dense linear subset of
Hilbert space $\mathcal{H}$, such that $q(\cdot,\psi)$ is linear and
$q(\varphi,\cdot)$ is skew-linear. If $q(\varphi,\psi) =
\overline{q(\psi,\varphi)}$ we say that $q$ is symmetric.

The map $q$ is called semi-bounded if $q(\varphi,\varphi)\geq
-M\|\varphi\|^2$ for some $M>0$. A semi-bounded quadratic form $q$ is
called closed if $Q(q)$ is complete under the norm
\[
   \|\psi\|_{+1} = \sqrt{q(\psi,\psi)+(M+1)\|\psi\|^2} .
\]

If $q$ is a closed semi-bounded quadratic form, then $q$ is the
quadratic form of a unique self-adjoint operator $A$, such that
$q(\varphi,\psi) = (A\varphi,\psi)$, for any $\varphi,\psi\in
D(A)$. Conversely, for any self-adjoint operator $A$ on $\mathcal{H}$,
there exists a corresponding quadratic form $q$. Then we denote
$Q(A)=Q(q)$. If $A$ is merely symmetric, let
$q(\varphi,\psi)=(A\varphi,\psi)$ for any $\varphi,\psi\in D(A)$. We
can complete $D(A)$ under the inner product
\[
   (\varphi,\psi)_{+1} = q(\varphi,\psi) + (\varphi,\psi)
\]
to obtain a Hilbert space $\mathcal{H}_{+1}\subset\mathcal{H}$. Then
$q$ extends to a closed quadratic form $\overline{q}$ on
$\mathcal{H}_{+1}$ (This is known as Friedrichs extension). We call
$\overline{q}$ the closure of $q$.

Let $A_1$, $A_2$ be two self-adjoint operators on Hilbert spaces
$\mathcal{H}_1$ and $\mathcal{H}_2$, with domains $D_1(A_1)$ and
$D_2(A_2)$. Then denote $A_1 \bigoplus A_2$ to be the operator $A$ on
$\mathcal{H}_1 \bigoplus \mathcal{H}_2$, with domain $D(A) = \{\langle
\varphi,\phi \rangle|\varphi\in D(A_1),\phi\in D(A_2)\}$ with
$A\langle \varphi,\psi \rangle = \langle A_1\varphi,A_2\phi \rangle$.
Now assume that $A_1$ and $A_2$ are nonnegative, and $\mathcal{H}_2
\subset \mathcal{H}_1$. We write $A_1 \leq A_2$ if and only if
\begin{enumerate}
\item $Q(A_1)\supset Q(A_2)$.
\item For any $\psi \in Q(A_2)$, $(A_1\psi,\psi) \leq (A_2\psi,\psi)$.
\end{enumerate}
~\\

\begin{lemma}
If $A_1$ and $A_2$ are two self-adjoint operators, $0\leq A_1\leq A_2$, then
\[\lambda_n(A_1)\leq\lambda_n(A_2),\]
where $\lambda_n(A)$ denotes the n-th eigenvalue of $A$, counted with multiplicity.
\end{lemma}

\medskip\medskip

\noindent
\textbf{Trace class}. Let $T:\mathcal{H}\rightarrow \mathcal{H}$ be a
compact operator, then $T^*T:\mathcal{H}\rightarrow \mathcal{H}$ is a
selfadjoint semidefinite compact operator, and hence it has a discrete
spectrum,
\[
   s_0(T)^2\geq s_1(T)^2\geq\cdots\geq s_k(T)^2\rightarrow 0 .
\]
We say that the nonnegative $s_j(T),j=1,2,\cdots$ are singular values
of $T$.\\ $T$ is said to be of trace class, if
\[
   \sum_{j=1}^\infty s_j(T)<\infty .
\]

\medskip\medskip
 
\begin{theorem}
Suppose that $T$ is of trace class on $\mathcal{H}$ and
$\{f_n\}_{n=1}^\infty$ is any orthonormal basis of $\mathcal{H}$, then
 \[\sum_{j=0}^\infty\langle Tf_j,f_j\rangle\]
converges absolutely to a limit that is independent of the choice of
orthonormal basis. The limit is called the trace of $T$, denoted by
$\mathrm{trace}(T)$.
\end{theorem}

\medskip\medskip
 
\begin{theorem}
Suppose that $T$ is of trace class on $\mathcal{H}$ and $B$ is a
bounded operator on $\mathcal{H}$, then $T B$ and $B T$ are of trace
class on $\mathcal{H}$, and
\[
   \mathrm{trace}(TB) = \mathrm{trace}(BT) .
\]
\end{theorem}

\medskip\medskip

\begin{theorem}
If $K$ is an operator of trace class on $L^2(\Omega)$, and
\[
   Kf(x)=\int_{\Omega}K(x,y)f(y)\mathrm{d}y
\]
then
\[
   \mathrm{trace}(K) = \int_{\Omega}K(x,x)\mathrm{d}x .
\]
\end{theorem}

\medskip\medskip

\noindent
Let $A(x,\epsilon D)$ be a semiclassical
pseudodifferential operator with scalar symbol $A(x,\xi)$, then
\[
   A(x,\epsilon D) f = \int_{\mathbb{R}^n} K_A(x,y) f(y) \mathrm{d}y ,
\]
with
 \[
   K_A(x,y) = \frac{1}{(2\pi \epsilon)^n}
     \int_{\mathbb{R}^n} A(x,\xi) e^{\frac{\ii}{\epsilon}(x-y)\cdot\xi}
                         \mathrm{d}\xi .
\]
If $A(x,\epsilon D)$ is of trace class, an immediate implication of
the theorem above is that
\[
   \mathrm{trace}({A(x,\epsilon D)})
         = \int_{\mathrm{R}^n} K_A(x,x) \mathrm{d}x
               = \frac{1}{(2\pi \epsilon)^n}
   \int_{\mathbb{R}^n} \int_{\mathbb{R}^n}
                            A(x,\xi)\mathrm{d}x\mathrm{d}\xi .
\]

{}
\end{document}